\documentclass{amsart}

\usepackage{amsfonts, amssymb} 
\usepackage{url}  
\usepackage{mathtools}
\usepackage{graphicx}
\usepackage[all]{xy}
\usepackage{mathrsfs}
\usepackage{color, marginnote}

\newtheorem{theorem}{Theorem}[section]
\newtheorem{lemma}[theorem]{Lemma}
\newtheorem{proposition}[theorem]{Proposition}
\newtheorem{corollary}[theorem]{Corollary}

\newtheorem{definition}[theorem]{Definition}

\theoremstyle{definition}
\newtheorem{remark}[theorem]{Remark}

\newcommand{\cL}{\mathcal{L}}
\newcommand{\cM}{\mathcal{M}}
\newcommand{\cA}{\mathcal{A}}
\newcommand{\cN}{\mathcal{N}}
\newcommand{\bR}{\mathbb{R}}
\newcommand{\R}{\mathbb{R}}

\newcommand{\bZ}{{\mathbb Z}}
\newcommand{\Z}{{\mathbb Z}}
\newcommand{\bS}{{\mathbb S}}
\newcommand{\bT}{{\mathbb T}}
\newcommand{\T}{{\mathbb T}}
\renewcommand{\a}{{\alpha}}
\renewcommand{\b}{{\beta}}
\newcommand{\g}{{\gamma}}

\newcommand{\m}{{\mu}}

\newcommand \fC {\mathfrak C}
\newcommand \calM {\mathfrak{M}}
\def\ie{\hbox{\it i.e.,\ }}
\newcommand \dpr {\partial}

\renewcommand{\epsilon}{\varepsilon}

\newcommand{\Hom}{\rm Hom}

\begin{document}

\title[On the set of asympt. homologies of orbits on invar. Lagr. graphs]{On the set of asymptotic homologies of orbits on invariant Lagrangian graphs}

\author{Rafael O. Ruggiero}
\address{Dept. de Matematica - PUC - Rio de Janeiro - Brazil}
\email{rafael.o.ruggiero@gmail.com }

\author{Alfonso Sorrentino}
\address{Department of Mathematics, University of Rome Tor Vergata, Via della ricerca scientifica 1, 00133 Rome (Italy).}
\email{sorrentino@mat.uniroma2.it}
\thanks{}



\date{\today}

\dedicatory{}

\begin{abstract} 
Given a smooth  Tonelli Hamiltonian on the torus $\mathbb{T}^{n}$ and a $C^{2}$ Lagrangian graph $W \subset T^{*}\mathbb{T}^{n}$ that is invariant under the Hamiltonian flow and contained within a Mañé supercritical energy level, we demonstrate the existence of a proper cone in the first real homology group $H_1(\T^n,\R)$ that contains the asymptotic homologies of the canonical projections of recurrent orbits in $W$. \\
Additionally, for invariant Lagrangian graphs on $\mathbb{T}^{3}$,  drawing on Franks' theory of the rotation set of homeomorphisms of $\mathbb{T}^{2}$ homotopic to the identity,
we show that under certain assumptions for an invariant Lagrangian graph on $\mathbb{T}^{3}$, if there exists a rational vector in homology contained in the set of asymptotic homologies of orbits on the Lagrangian graph, then the graph contains  a Mather measure supported on a periodic orbit. This result generalizes a well-known fact for Lagrangian graphs on $\mathbb{T}^{2}$.
Finally, we exploit these results for three dimensional tori to give a partial answer to a conjecture by Carneiro-Ruggiero about the non-existence of Hedlund Lagrangian tori at supercritical energy levels.
\end{abstract}

\maketitle

\section{Introduction} 
From the seminal works of Morse and Hedlund on globally minimizing geodesics in the universal cover of compact surfaces of non-negative genus, it is well established that globally minimizing geodesics behave analogously to geodesics in the hyperbolic plane when the surface has genus greater than one, and similarly to straight lines in the Euclidean plane when the surface is a torus. In the latter case, it is natural to associate each globally minimizing geodesic with an asymptotic homology class, which corresponds to a vector tangent to a straight line that ``shadows" the geodesic in the Euclidean plane (the so-called {\it Schwartzman asymptotic cycle}). This elegant geometric feature of globally minimizing geodesics, however, is not generally observed in higher-dimensional tori, where the asymptotic homology of global minimizers may fail to exist.\\

One of the key contributions of John Mather, foundational to what is now known as {\it Aubry-Mather theory}, is the observation that it is still possible to associate an asymptotic homology class with certain special orbits and invariant probability measures for so-called {\it Tonelli} Hamiltonian system on the cotangent bundle of any $n$-dimensional torus (or, more generally, any closed manifold).\\

Recall that a function $H:\,T^*\mathbb{T}^n\longrightarrow \mathbb{R}$ is called a {\em Tonelli Hamiltonian} if:
\begin{itemize}
    \item[i)]  $H$ is of class $C^2(T^*\T^n)$;
    \item[ii)]  $H$ is strictly convex in each fiber in the $C^2$ sense, {\it i.e.}, the second partial vertical derivative $\frac{\partial^2 H}{\partial p^2}(x,p)$ is positive definite as a quadratic form, for any $(x,p)\in T^*\mathbb{T}^n$;
    \item[iii)] $H$ is superlinear in each fiber, {\it i.e.},
    $$\lim_{\|p\|_x\rightarrow +\infty} \frac{H(x,p)}{\|p\|_x} = + \infty\,.$$
\end{itemize}

The special orbits or measures central to Aubry-Mather theory are obtained by minimizing the so-called {\it Lagrangian action} (see Appendix \ref{sec2.4} for a more precise description) and in a sense, generalize orbits or invariant measures on invariant Lagrangian graphs. Indeed, it is not difficult to show that orbits or invariant measures on an invariant Lagrangian graph are action-minimizing. 

This observation allowed Mather to construct a rich variety of compact invariant subsets of the system (known as, {\it Aubry-Mather sets}), which, while not generally submanifolds or regular in a classical sense, can be viewed as generalizations of invariant Lagrangian graphs. Besides their dynamical significance, these objects also play an important role in other contexts, such as PDEs (e.g., the Hamilton-Jacobi equation and weak KAM theory), symplectic geometry, and beyond (see for example \cite{Sorrentinobook} for more details).\\

Understanding the relationship between these Aubry-Mather sets and the asymptotic homology classes of action-minimizing orbits is a fundamental question, which lies at the heart of several interesting results and open problems in the field.
When the Hamiltonian flow of a Tonelli Hamiltonian is integrable in the sense of the Arnol'd-Liouville theorem, Aubry-Mather sets correspond to invariant Lagrangian graphs (which are diffeomorphic to tori) that foliate the phase space. In particular, all orbits on a given leaf ({\it i.e.}, an invariant Lagrangian graph) share the same asymptotic homology class. A similar phenomenon occurs under weaker notions of integrability (see for example \cite{Arnaud, BS, MS}).\\

This setting naturally raises many questions about the asymptotic homologies of orbits lying on invariant Lagrangian graphs that do not necessarily belong to a foliation. Specifically: {\it Is there a reasonable description of the asymptotic homology of orbits or Mather measures in a Lagrangian invariant graph?}\\

In this article, we address this question by proving the following:

\begin{theorem} \label{main}
Let $H: T^{*}\mathbb{T}^n \longrightarrow \mathbb{R}$ be a Tonelli Hamiltonian on the torus, and let $W \subset T^{*}\mathbb{T}^{n}$ be a $C^{2}$ Lagrangian graph invariant under the Hamiltonian flow at a supercritical energy level $E > c_{0}(L)$, where $L$ is the corresponding Lagrangian and $c_{0}(L)$ is the strict Ma\~{n}\'{e} critical value (see Definition \ref{strictcL}). Then, there exists a proper cone $C$ in $H_{1}(\mathbb{T}^{n}, \mathbb{R})$, with respect to the stable norm in homology, such that the asymptotic homologies generated by quasi-orbits of the canonical projection of the Hamiltonian flow restricted to $W$ are contained within $C$.
\end{theorem}

For the definition of quasi-orbit we refer to Definition \ref{quasiorbit}.

Theorem \ref{main} reveals a dichotomy in the geometry of asymptotic homology for Lagrangian invariant graphs. If the Lagrangian graph is supercritical in the sense of Theorem \ref{main}, then the asymptotic homologies are confined to a positive cone in homology, whose size, interestingly, is uniform across the category of Riemannian metrics and Mechanical Hamiltonian systems (see Proposition \ref{compactcase} and Remark \ref{remmechanical}). Conversely, if the Lagrangian graph lies at Ma\~n\'e's strict critical value, no clear control exists over the geometry of the asymptotic homology; we describe an example in Section \ref{maneexample}.\\

We also explore the following question, which is related to the structure of Mather's minimal average actions (the so-called Mather's $\alpha$ and $\beta$ functions, see Appendix \ref{sec2.4} for a precise definition):\\

\noindent {\sc Question:} {\it If the asymptotic homologies of the orbits of the Lagrangian graph are all irrational, does the cone $C$ have an empty interior?}\\

We investigate this problem for the case where the dimension of $W$ is $3$, drawing on Franks' theory of the rotation set of homeomorphisms of $\mathbb{T}^{2}$ homotopic to the identity. We prove the following result, discussed in Section \ref{secHedlund}:

\begin{theorem} \label{main2}
Let $H: T^{*}\mathbb{T}^n \longrightarrow \mathbb{R}$ be a Tonelli Hamiltonian on the torus, and let $W \subset T^{*}\mathbb{T}^{n}$ be a $C^{2}$ Lagrangian graph invariant under the Hamiltonian flow at a supercritical energy level $E > c_{0}(L)$, where $L$ is the corresponding Lagrangian and $c_0(L)$ denotes Ma\~{n}\'{e}'s strict critical value (see Definition \ref{strictcL}). Let $\psi_{t} :  \mathbb{T}^n \longrightarrow \mathbb{T}^n$ be the canonical projection of the Hamiltonian flow restricted to $W$. Then $\psi_{t}$ is conjugate to the suspension of a homeomorphism $f : \mathbb{T}^{n-1} \longrightarrow \mathbb{T}^{n-1}$. Moreover, if $n=3$ and $f$ is homotopic to the identity, the absence of periodic orbits implies that the minimal proper cone containing all asymptotic homologies of quasi-orbits of $\psi_{t}$ is contained in a codimension-1 plane.
\end{theorem}

The same ideas involved in the proof of Theorem \ref{main2} allow us to show that, for three dimensional tori, if the cone generated by the asymptotic homologies of quasi-orbits is an open set containing a rational vector, then the Lagrangian graph contains a Mather measure supported on a periodic orbit, provided that the associated suspension is homotopic to the identity. This result generalizes a well-known fact for Lagrangian graphs on $\mathbb{T}^{2}$.\\

Moreover, we exploit Theorem \ref{main2}  to show the non-existence of so-called {\it Hedlund Lagrangian graphs} in three-dimensional tori,
under the assumptions of Theorem \ref{main2}. Hedlund, in fact,  constructed in \cite{kn:Hedlund} examples of Riemannian metrics on the $n$-torus such that the only Mather measures ({\it i.e.}, action-minimizing invariant measures) are supported on $n$ closed geodesics whose homologies generate the first homology group. An interesting question arising from Hedlund's examples is whether there exists a Lagrangian invariant submanifold for the geodesic flow of a Hedlund metric in the torus such that the only Mather measures are supported on $n$ closed orbits whose homologies generate the first homology group. We shall call such a Lagrangian submanifold {\it Hedlund Lagrangian submanifold} after \cite{kn:CarRug} (see also Definition \ref{Hedlundtorus}).\\

 The existence of Hedlund Lagrangian submanifold is an intriguing question. Carneiro and Ruggiero studied this problem in \cite{kn:CarRug}, proving that a Hedlund Lagrangian submanifold must be a graph (called a {\it Hedlund Lagrangian graphs} in \cite{kn:CarRug}), concluding their the generic non-existence  at supercritical energy levels (see also Proposition \ref{Hedlund-tori-nonex1} and discussion thereafter). They conjecture that Tonelli Hamiltonians on the $n$-dimensional torus do not admit Lagrangian invariant tori at supercritical energy levels where the set of Mather measures consists of exactly $n$ measures, and their asymptotic homologies generate the full homology group $H_1(\mathbb{T}^n,\mathbb{R})$.

In this article, we prove Carneiro-Ruggiero conjecture in dimension 3 under certain assumptions for the Hedlund Lagrangian submanifold (see Proposition \ref{propmain3}).\\

\subsection*{Organization of the paper}
The article is organized as follows:
\begin{itemize}
\item In Section \ref{secprelim} we recall  some preliminary material that will be relevant in the following. More specifically: Finsler metrics and their properties (Section \ref{secFinsler}), 
Stable norm and asymptotic homology (Section \ref{secStablenorm}), Finsler geometry of supercritical energy levels of Tonelli Hamiltonians (Section \ref{secFinHam}), Hamilton-Jacobi equation (Section \ref{secHJ}), and some facts on foliations defined by closed one forms (Section \ref{secFoliation}).
\item In Section \ref{secProofmainthm} we prove Theorem  \ref{main}. In particular, in Section \ref{secfinslertoRiem} we show that the proof of Theorem \ref{main} follows from the proof of the same statement for Riemannian geodesic flows, while in Section \ref{sechamtoriem} we discuss how to relate the Hamiltonian case to the Riemannian geodesic flow case. The proof of Theorem \ref{main} will then follow from the one of Theorem \ref{cone-fib} in Section \ref{secconehom}.
{In Section \ref{maneexample} we discuss an example to show that Theorem \ref{main} might fail on the energy level corresponding to Ma\~n\'e's strict critical value.}
\item In Section \ref{secdim3} we describe  how to use Theorem \ref{main}  to prove  Theorem \ref{main2}. In particular, the proof will be based on the relationship between the rotation set of a homeomorphism of $\mathbb{T}^{2}$ that is homotopic to the identity, and the asymptotic homologies of quasi-orbits of its suspension flow in $\mathbb{T}^{3}$. The relevant backgroud material will be recalled in Section \ref{sec5-1} and the proof of the theorem will be provided in Section \ref{sec5-2}. Moreover, in Section \ref{sec5-3} we describe an algebraic criterion that implies the uniqueness of the asymptotic homologies of Mather measures in the graph (see Lemma \ref{totally-irrational}). 
\item In Section \ref{secHedlund} we show the non-existence of Hedlund Lagrangian graphs in three dimensional tori in the same setting  as in Theorem \ref{main2} (see Proposition \ref{propmain3}).
\item In Appendix \ref{sec2.4} we provide a short introduction to Aubry-Mather theory for Lagrangian systems.
\item in Appendix \ref{Appcorrigendum} we include a corrigendum  to the proof of  \cite[ Lemma 3.3]{kn:CarRug}.\\
\end{itemize}

\subsection*{Acknowledgements}
RR was partially supported by  PRONEX de Geometria, FAPERJ, E-16/2014 INCT/FAPERJ, Projeto Universal CNPq 430154/2018-6, and Mathamsud CAPES/CDEFI n.  88881.878892/2023-01. 
AS  acknowledges the support of the Italian Ministry of University and Research’s  PRIN 2022 grant ``{\it Stability in Hamiltonian dynamics and beyond}’', as well as the Department of Excellence grant MatMod@TOV (2023-27) awarded to the Department of Mathematics of University of Rome Tor Vergata.  AS is a member of the INdAM research group GNAMPA and of the UMI group DinAmicI.\\

\section{Preliminaries}\label{secprelim}

For the reader's sake, in this section we would like to recall some preliminary material that will be relevant in the subsequent discussion and the proof of Theorem \ref{main}.

\subsection{Finsler metrics} \label{secFinsler}

For more details and proofs of the results contained in this subsection, we refer the reader to \cite{kn:BCS}. 

Let $M$ be a $n$-dimensional $C^{\infty}$ manifold, let $T_{x}M$ be
the tangent space at $x\in M$, and let $TM$ denote its tangent bundle.
In local coordinates, an element of $T_{x}M$ can be expressed as a
pair $(x,v)$, where $v$ is a tangent vectorat $p$. 
Let $\pi:TM \longrightarrow M$, $\pi(x,v) := x$ denote the canonical projection. 
Moreover, we denote by $\Sigma_0 := \{(x,v) \in TM; \, v = 0 \} $ the zero-section.

\begin{definition}
A $C^{k}$ ($k\geq 2$) {\em Finsler structure} (or {\em metric}) on $M$ is a function
$F:TM \rightarrow [0,+ \infty)$ satisfying the following properties: 
\begin{itemize}
\item[(i)] $F$ is $C^k$ on $TM\setminus \Sigma_0$; 
\item[(ii)] $F$ is positively homogeneous of degree one in $v$, {\it i.e.},
$$ F(x,\lambda v)=\lambda\, F(x,v) \qquad \forall \, \lambda >0, \quad \forall\,(x,v) \in TM;$$
\item[(iii)] The Hessian matrix of $F^2$, namely
$$g_{ij} := \frac{1}{2}\frac{\partial ^2}{\partial v^i \partial v^j}F^2  \qquad i,j\in \{1,\ldots, n\}$$
is positive definite on $TM\setminus \Sigma_0$.
\end{itemize}
A {\em Finsler manifold} is a pair $(M,F)$ consisting of a $C^{\infty}$ manifold $M$ and a
$C^{\infty}$ Finsler structure $F$ on $M$.
\end{definition}

\begin{remark}
{\bf (i)}  Any Riemannian metric is of course a Finsler metric, where the function $F^{2}(x,v)$ is a quadratic polynomial in the fiber variables $v$.\\
{\bf (ii)} If $M$ is compact, it is not difficult to show (see \cite{kn:BCS} ) that given any two Finsler metrics $F_{1}$, $F_{2}$, 
there exists a constant $C>0$ such that 
$$ \frac{1}{C} F_{1}(x,v) \leq F_{2}(x,v) \leq C F_{1}(x,v) \qquad \forall\, (x,v)\in TM.$$
Therefore, any two Finsler metrics on a compact manifold $M$ are equivalent. 
\end{remark}

The strict convexity of the Hessian of $F^{2}$ allows one to apply the Lagrangian formalism of calculus of variations and introduce the corresponding geodesic flow $\phi_{t} : TM \longrightarrow TM$ given by  the Euler-Lagrange flow associated to $(M,F)$.  \\
We call any canonical projection of an
orbit of this geodesic flow a  {\it Finsler geodesic}. 
For a non-vanishing vector $v\in T_{x}M$, we shall denote by
$\gamma _{(x,v)}(t)$ the geodesic with initial conditions $\gamma
_{(x,v)}(0)=x$ and $\dot{\gamma}_{(x,v)}(0)=v$. \\

As a Lagrangian flow, the geodesic flow is homogeneous in the tangent bundle, therefore the Lagrangian flow restricted to an energy level  is a reparametrization of the Lagrangian flow on the unit tangent bundle, that we shall denote by 
(the unit tangent bundle depends on $F$, we shall omit this dependence in the notation):

$$T^{1}M := \{ (x,v) \in TM:\; F(x,v) =1\}.$$ 
Hereafter, we shall assume that all geodesics have unit speed.  \\

The length of an absolutely continuous curve  $\xi:[a,b]\rightarrow M$ with respect to the Finsler 
structure $(M,F)$ is given by 
$$ \ell_F(\xi):= \int _a^b F(\xi(t), \dot{\xi}(t))\, dt.$$
$\ell_F$ gives rise to what is called a \textit{Finsler distance}, namely, a   function  
\begin{eqnarray*}
d_F:M\times M &\longrightarrow&
[0,\infty)\\
(x,y) &\longmapsto& d_F(x,y) := \inf _\xi {\ell_F(\xi)},
\end{eqnarray*}
where the infimum is taken over all absolutely continuous curves
$\xi:[0,1] \rightarrow M$ with $\xi(0)=x$ and $\xi(1)=y$. 
However, {differently from} the Riemannian case, $d_F(x,y)$ might 
be different from $d_F(y,x)$ and therefore, $d_{F}$ might not be a distance in the usual sense. \\

A Finsler manifold $(M,F)$ induces naturally a Finsler structure
in the universal covering $\tilde{M}$ of $M$, just by pulling back
the Finsler structure $F$ to the tangent space of $\tilde{M}$ by the
covering map. Let us denote by $(\tilde{M},\tilde{F})$ this Finsler
manifold. Then $(\tilde{M},\tilde{F})$ and $(M,F)$ are locally isometric.\\

A geodesic $\sigma: [a,b] \longrightarrow \tilde{M}$ is called
\textit{forward minimizing}, or simply minimizing, if
$\ell_{\tilde{F}}(\sigma) \leq\ell_{\tilde{F}}(\xi)$ for all  absolutely continuous curves
$\xi:[a,b]\longrightarrow \tilde{M}$ such that $\xi(a)=\sigma(a) $,
$\xi(b)=\sigma(b)$ (this implies that $\sigma: [s,t] \longrightarrow
\tilde{M}$ is also minimizing for every $a\leq s \leq t \leq b$). 
Notice that, in general, a minimizing geodesic $\sigma $ might fail
to be minimizing if one reverses its orientation.\\

\subsection{Stable norm and asymptotic homology}\label{secStablenorm}

Given a $C^{\infty}$ manifold $M$, we shall denote by $H_{1}(M,\mathbb{Z})$ and $H_{1}(M,\mathbb{R})$
the first homology groups of $M$ with, respectively, integer and real; similarly, for the first cohomology groups $H^{1}(M,\mathbb{Z})$ and $H^{1}(M,\mathbb{R})$. 

{In this section we introduce {Federer-Gromov's} {\it stable norm} on $H_1(M;\bR)$ and state some of its properties. We refer the reader to \cite{BBI, GLP} for a more 
detailed presentation.

\begin{definition}
Let $(M,F)$ be a $C^{\infty}$ Finsler structure of $M$ and let $\ell_{F}(\gamma)$ denote the length of an absolutely continuous curve $\gamma$ with respect to the metric $F$. Given a closed curve $\beta : \bS^{1} \longrightarrow M$, let $[\beta] \in  H_{1}(M,\mathbb{Z})$ denote its homology class.\\
 For every $h \in H_{1}(M,\mathbb{R})$, we define the {\em stable norm} $\|h\|_{s}$ as:
 $$ \| h \|_{s} :=  \inf \left\{ \sum_{i=1}^{m} r_{i}\ell_{F}(\gamma_{i}): \; \mbox{$\gamma_{i}$ is a closed curve   and $h = \sum_{i=1}^{m}r_{i}[\gamma_{i}]$, $r_i \geq 0$}
\right \}. $$
\end{definition}

\begin{remark}
{\bf (i)} The stable norm clearly depends on the metric.\\
{\bf (ii)} If $M=\mathbb{T}^{n}$, it is a norm in $H_1(\mathbb T^n,\mathbb R)\simeq \mathbb{R}^{n}$ and therefore, the stable norms associated to two Finsler metrics in $\mathbb{T}^{n}$ are equivalent. \\
{\bf (iii)} We also have that the stable norm of an integer homology class in $(M,F)$ is the length of a closed geodesic that minimizes length in the family of closed curves in the class. 
\end{remark}

\medskip

Let us now introduce what we shall call a \textit{cone} in $H_{1}(M,\mathbb{R})$ relative to the stable norm associated to $(M,F)$. 

\begin{definition}
Le  $h \in H_{1}(M,\mathbb{R})$ be a non-trivial class and let $\langle h \rangle$ denote the line generated by $h$ in $H_{1}(M,\mathbb{R})$. Given a decomposition as direct sum  $H_{1}(M, \mathbb{R}) = \langle h \rangle \oplus \mathcal{G}$, and $A \geq 0$, we define the {\em cone of slope $A$ with axis $h$} as 
$$ \fC_{A}(h) := \bigcup_{t\geq0}\left\{(th, w) \in \langle h \rangle \times \mathcal{G}: \;  \frac{\| w \|_{s}}{\| t h \|_{s}} \leq A \right\}. $$
\end{definition}

\begin{remark}
{\bf (i)} The direct sum decomposition $H_{1}(M, \mathbb{R}) = \langle h \rangle \oplus \mathcal{G}$ is not unique, so the cone $\fC_{A}(h)$ might depend on the chosen direct sum. Notice, however, that  a cone is always properly contained in a half semi-space of $H_{1}(M, \mathbb{R})$, independently of the chosen direct sum decoposition. 
 Clearly, different metrics may define different  cones $\fC_{A}(h)$ with respect to a given direct sum decomposition of $H_{1}(M, \mathbb{R})$.\\
{\bf (ii)} If $A=0$, then we just get $\fC_{A}(h)$ is the half-line $\{ th:\; t\geq 0\}$.  \\
\end{remark}

We finish this subsection by recalling a slight extension of the Schwartzmann's notion of \textit{asymptotic homology}  (see \cite{kn:Schwartz}), also called \textit{homological rotation vector},
given by Burago-Ivanov-Kleiner in \cite{kn:BurIvaKle} for Riemannian metrics.

\begin{definition}\label{defhomvector}
Let $(M,F)$ be a $C^{\infty}$ complete Finsler structure on $M$ and consider $\xi:[0, \infty) \longrightarrow M$ a uniformly continuous curve.
For every $t>0$, denote by $\gamma_{t}: [0,1] \longrightarrow M$ be {shortest curve} joining $\xi(t)$ to $\xi(0)$. Let
$\Xi_{t}$ be the closed curve formed by the concatenation of $\xi:[0,t] \longrightarrow M$ and $\gamma_{t}$ and denote by
$[\Xi_{t} ]\in H_1(M,\bZ)$  the integer homology class of the curve $\Xi_{t}$. Then, if the limit 
$$  [\xi]:= \lim_{t \rightarrow \infty}\frac{[\Xi_t]}{t}  $$ 
exists, then it is called the {\em asymptotic homology class} or {\rm homological rotation vector} of $\xi$. \end{definition}

\begin{remark}
{\bf (i)} It is easy to show that the asymptotic homology class is independent of the choice of the metric $F$, 
see for instance \cite{kn:CarRug}.\\
{{\bf (ii)} Examples of curves with well defined asymptotic homology classes arise in the study of global action minimizers of Tonelli Lagrangians (see next subsection for the definition) and in what is called {\it Aubry-Mather theory}; we refer the reader to  Appendix A and \cite{Gonzalonotes, Mane2, Ma, Sorrentinobook} for more details.} 
 \end{remark}

\begin{definition}
We say that the curve $\xi: [0, \infty) \longrightarrow M$, $t \in \mathbb{R}$,  has {\em totally irrational} asymptotic homology  if there exists a basis $\{h_{i}\}_{i=1}^m$ of $H_{1}(M,\mathbb{Z})$ such that the coefficients of the asymptotic homology class  in this basis, namely $[\xi] = (\alpha_1,\ldots, \alpha_m) \in \bR^m$ form a set of rationally independent number (this means that if there exist rational numbers $\{a_{1},a_{2},..,a_{m}\}$ such that  $\sum_{i=1}^{m}a_{i}\alpha_{i} =0$, then $a_{i}=0$ for every $i=1,\ldots, m$). 
\end{definition}

\begin{remark}
{\bf (i)} Notice that the property of having a totally irrational homology does not depend on the choice of a basis for $H_{1}(M, \mathbb{Z})$, since every  basis  can be obtained from a given one applying a linear transformation with rational entries. \\
{\bf (ii)} Clearly, a curve with totally irrational asymptotic homology cannot be a closed curve. 
\end{remark}

 \medskip
\subsection{Finsler geometry of supercritical energy levels of Tonelli Hamiltonians}\label{secFinHam}

Let us first briefly recall the definition of the \textit{Ma\~n\'e's critical value}
(see, for instance, \cite{kn:CIPP98, Mane2} for more details). \\

We consider a Tonelli Lagrangian $L: TM\longrightarrow \bR$ and the associated Tonelli Hamiltonian $H:T^*M\longrightarrow \bR$.\\

Given any two points $(x,y) \in M$ and $T>0$, denote by $\mathcal{C}_T(x,y)$ the set of absolutely continous curves $\gamma :[0,T] \rightarrow M$, with $\gamma (0)=x$ and $\gamma (T)=y$. We define, for each $k \in \mathbb{R}$, the \textit{action potential} $\Phi _k: M \times M \rightarrow \mathbb{R}$ by 
$$ \Phi _k(x,y)= \inf \{ A_{L+k}(\gamma ): \; \gamma \in \cup _{T>0}\mathcal{C}_T(x,y) \} , $$
where 
$$ A_{L+k}(\gamma )=\int _a^b \left(L(\gamma (t), \gamma '(t)) + k\right)\, dt $$
denotes the action of $L+k$ on $\gamma $.

\begin{definition}
The {\em Ma\~n\'e's critical value} of $L$ is the real number $c(L)$ defined as the infimum of $k \in \mathbb{R}$ such that $\Phi _k(x,x) > - \infty $ for some $x \in M$.
\end{definition}

\medskip

\begin{remark}
It is easy to check that the above infimum is not $-\infty$. In fact, since $L$ is superlinear, there exists $k\in \bR$ such that $L+k\geq 0$ on $TM$ and therefore $\Phi_k(x,x)\geq 0$ for every $x\in M$.\\
\end{remark}

Let $E:TM \rightarrow \mathbb{R}$ be the energy function of $L$, that is, 
$$E(x,v)=\frac{\partial L}{\partial v}(x,v)\cdot v - L(x,v) ,$$
for all $(x,v) \in TM$ (in other words, the Hamiltonian pulled back on $TM$ via the Legendre transform \eqref{FenchelLeg}). Let $\Phi^L_t:TM \rightarrow TM$ be the Euler-Lagrange flow of $L$.\\

Let us denote by $\tilde M$ the abelian covering of $M$ and let $\Pi: \tilde{M} \longrightarrow M$ be the covering map. We can lift the Lagrangian $L$ to $T\tilde{M}$ by defining $\tilde{L}(x,v) = (L \circ \Pi)(x,v)$. 

\begin{definition}\label{strictcL}
The {\em strict Ma\~n\'e's critical value} of $L$, $c_{0}(L)$, is given by 
$c_{0}(L) := c(\tilde{L})$. 
\end{definition}

The following proposition is contained in \cite[Corollary 2]{kn:CIPP98}, shows the relevance of the geodesic flow in the study of Tonelli Hamiltonian flows.

\begin{proposition}\label{propfinslerHam}
If $k>c_{0}(L)$, then the flow $\varphi _t|_{E^{-1}(k)}$ is a reparametrization of the geodesic flow on the unit tangent bundle of an appropriately chosen Finsler metric on $\tilde M$. 
\end{proposition}

\bigskip

\subsection{Hamilton-Jacobi equation}\label{secHJ}

Let $\omega=-d\alpha$ denote the symplectic form on $T^*M$, where $\alpha$ -- which in local coordinates reads $\sum_{i=1}^{} p_i dx_i$, where $n=\dim M$ -- is the {\it Liouville form} (or {\it tautological form}). \\
Let us recall some definitions.

\begin{definition}
A $C^{k}$ submanifold $S \subset T^{*}M$ is called {\em Lagrangian} if it has the same dimension of $M$ and for every $\theta \in S$, the restriction of the canonical $2$-form $\omega$ to the tangent space $T_{\theta}S$ vanishes. 
 \end{definition}

\begin{remark}
{\bf (i)} One can check that the property of $S \subset T^{*}M$ being Lagrangian is equivalent to the condition that the restriction of the canonical $1$-form $\alpha$ to $S$ is a closed $1$-form. \\
{{\bf (ii)} It is easy to check that if $S \subset T^{*}M$ is Lagrangian and invariant under the Hamiltonian flow of $H$, then $H$ is constant on each connected component of $H$ (hence, its connected components lie in energy levels of the Hamiltonian).
}
\end{remark}

We say that a Lagrangian submanifold $W \subset T^{*}M$ is a \textit{Lagrangian graph} if the restriction of the canonical projection $\pi: T^*M \longrightarrow M$, $\pi(x,p)=x$,
 to $W$ is a diffeomorphim. If $W$ is a Lagrangian graph that is invariant by a Hamiltonian flow $\Phi^H_{t}: T^{*}M \longrightarrow T^{*}M$, then it is clear that the canonical projection of $W$ endows $M$ with a flow $\psi_{t} : M\longrightarrow M$ that is everywhere tangent to the canonical projection of the orbits of $\Phi^H_{t}$. In the case of geodesic Hamiltonians associated to a Riemannian metric, the projected orbits are just geodesics for the Riemannian metric. \\

A Lagrangian graph then endows the manifold $M$ with a closed 1-form that is invariant by the flow $\psi_{t}$. This special property of Lagrangian graphs connects such graphs with the so-called \textit{Hamilton-Jacobi} equation. \\

The Hamilton-Jacobi equation associated to the Hamiltonian $H: T^{*}M \longrightarrow \mathbb{R}$ is given by
$$ H(q, d_{q}f) =c $$
where $f: M\longrightarrow \mathbb{R}$ is a smooth scalar function and $c$ is a real constant. If such a function exists it is called a solution of the Hamilton-Jacobi equation. {It is easy to check that if $M$ is compact, there can be at most one value of $c$ for which smooth solutions exist}.
However, smooth solutions of the Hamilton-Jacobi equation do not exist in general and one looks for less regular solutions; the study of generalized solutions called {\it viscosity solutions} (related to what is known as {\it weak KAM theory}) has been object of intense research in partial differential equations and more recently in calculus of variations and Hamiltonian dynamics; see for instance \cite{Fathi}.

\subsection{Foliations defined by closed one forms} \label{secFoliation}
Observe that if $W$ is a Lagrangian graph in $T^*M$, then it corresponds to the graph of a closed $1$-form $\sigma_W$ on $M$ (see for instance \cite[Proposition A.1.2]{Sorrentinobook}). Moreover, if $W$ is invariant by the Hamiltonian flow of $H$, then it is contained in a energy level 
$H=c$  (see for instance \cite[Proposition A.1.5]{Sorrentinobook}) and it satisfies therefore the equation $H(q,\sigma_{W,q}) = c$. 
Since every closed 1-form is locally exact, for every given $q \in M$ there exists an open neighborhood $U_{q} \subset M$ such that $\sigma_W = d f$, for some scalar function $f: U_{q} \longrightarrow \mathbb{R}$. The function $f$ is a local primitive of the 1-form $\sigma_W$, the Hamilton-Jacobi equation has therefore a local solution since 
$$ H(x, \sigma_W) = H(x, d_{x}f)= c$$ 
for $x \in U_{q}$. Since every two different primitives of $\sigma_W$ differ by a constant, the level sets of the primitives are well defined in $M$, and if $\sigma_W$ is nonsingular at every point this gives rise to a foliation $\mathcal{F}$ tangent to the kernel of $\sigma_W$.

 The theory of foliations defined by closed 1-forms is a classical subject in foliation theory, and the following results characterize precisely the foliation.
(see for instance \cite{kn:Molino}).
\begin{theorem} \label{fol}
Let $M$ be a compact, $n$-dimensional manifold endowed with a nonsingular closed 1-form $\sigma$. Then, a foliation 
$\mathcal{F}$ that is always tangent to the kernel of $\sigma$ satisfies one of the following properties:
\begin{enumerate}
\item[(i)] Every leaf is compact and diffeomorphic to some manifold $N$ whose dimension is $n-1$. 
\item[(ii)] Every leaf is dense. 
\end{enumerate}
\end{theorem}

The next wellknow results are due to Tischler \cite{kn:Tischler}.

\begin{theorem} \label{tischler}
Let $M$ be a $C^{\infty}$ compact manifold admitting a $C^{2}$ closed 1-form that does not vanish. Then $M$ fibers over $\bS^{1}$. Moreover, the closed 1-form can be arbitrarily approximated by closed 1-forms defining foliations by closed leaves. 
\end{theorem}

\begin{corollary} \label{fibred-torus}
Let $M = \mathbb{T}^{n}$ and let $\alpha$ be a $C^{2}$ closed 1-form defined on $\mathbb{T}^{n}$ without singularities. Then, given $\epsilon >0$, there exists a closed 1-form $\alpha_{\epsilon}$ without singularities such that $\alpha_{\epsilon}$ is $\epsilon$-close to $\alpha$ in the $C^{2}$ topology, and the foliation $\mathcal{F}$ defined by $\alpha$ is $\epsilon-$ close, in the $C^2$ topology, to a fibration by $n-1$ dimensional tori whose fibers are tangent to the kernel of $\alpha_{\epsilon}$.
\end{corollary}

\begin{proof} 
The proof follows at once from Theorem \ref{tischler} 
and the fact that the fundamental group of each leaf is a subgroup of the fundamental group of $\mathbb{T}^{n}$. Therefore, the fundamental group of the fibers is an abelian subgroup of $\mathbb{Z}^{n}$ of less rank, and the codimension one of the leaves implies that it must be isomorphic to $\mathbb{Z}^{n-1}$. In particular, the leaves are $(n-1)$ dimensional tori. 
\end{proof}

\bigskip

\section{Proof of Theorem \ref{main}} \label{secProofmainthm}

The purpose of the section is to show Theorem \ref{main}, the main ideas of the proof are presented in three subsections and the complete proof is given at the end of subsection 3.3. 

\subsection{Reduction of the Finsler case to the Riemannian case} \label{secfinslertoRiem}

The goal of the subsection is to show that the proof of Theorem \ref{main} follows from the proof of the same statement for Riemannian geodesic flows. Let us start with the following remark going back to the works of Gluck \cite{kn:Gluck} and Sullivan \cite{kn:Sullivan} about geodesible flows.

\begin{lemma} \label{Geodesible}
Let $(M,\bar{g})$ be a compact $C^{\infty}$ Riemannian manifold, and let $\mathcal{F}$ be a $C^{2}$ foliation whose tangent subbundle is the kernel of a nondegenerate closed 1-form $\beta$ of class $C^{2}$. Let $X$ be a nonsingular vector field such that $\beta_{x}(X(x) ) \neq 0$ for every $x \in M$. Then, there exists a Riemannian metric $g$ on $M$ such that the integral orbits of $X$ are geodesics of $(M,g)$, and the leaves of $\mathcal{F}$ are locally equidistant with respect to this metric. In particular, the integral flow of the normalized vector field $\bar{X}(x) = \frac{X(x)}{\parallel X(x) \parallel}$ preserves the foliation $\mathcal{F}$. 
\end{lemma}

\begin{proof}
We sketch the proof of Lemma \ref{Geodesible} for the sake of completeness (see \cite{kn:Sullivan} for more details). The metric $g$ can be obtained as follows. Given $x \in M$, let $\mathcal{F}(x)$ the leaf  containing 
$x$ and let $T_{x}\mathcal{F}(x)=\ker(\beta_{x})$  be the tangent space at $x$ of $\mathcal{F}({x})$. We have that the tangent space at $x$ equals te direct sum  
$$T_{x}M= \langle \lambda X(x) \rangle \oplus T_{x}\mathcal{F}(x),$$ 
where $\lambda \in \mathbb{R}$, and $\langle \lambda X(x) \rangle$ is the subspace generated 
by $X(x)$. 

Let $ P^{H}: T_{x}M\longrightarrow \langle \lambda X(x) \rangle$ be the projection onto the subspace $\langle \lambda X(x) \rangle$ and let $P^{V}: T_{x}M\longrightarrow \ker(\beta_{x})$ be the projection onto the kernel of $\beta_{x}$. The metric $g$ can be obtained as a product metric 
$$ g_{x} =P^{H}_{*}\bar{g}_{x} + P^{V}_{*}\bar{g}_{x} $$
This Riemannian metric is clearly $C^{2}$ if the flow $\psi_{t}$ and the foliation $\mathcal{F}$ are of class $C^{2}$. \\
The vector field $X$ and the tangent spaces of the leaves of $\mathcal{F}$ are perpendicular with respect to $g$. Then the orbits of $\psi_{t}$ are geodesics of $g$ by elementary variational arguments (see \cite{kn:RugGeoDed} for instance). In particular, the parameter $t$ of the integral flow $\bar{\psi}_{t}$ of the normalized vector field $\bar{X}(x)$ is locally the distance between the leaves of $\mathcal{F}$, therefore, 
$$ \bar{\psi}_{t}(\mathcal{F}(x)) = \mathcal{F}(\bar{\psi}_{t}(x))$$ 
for every $x \in M$. 

\end{proof}

\bigskip

\subsection{From Tonelli Hamiltonian to a Riemannian geodesic flow} \label{sechamtoriem}

Proposition \ref{propfinslerHam} and Lemma \ref{Geodesible} imply that it suffices to show Theorem \ref{main} for geodesic flows of Riemannian metrics. 

\begin{proposition} \label{Riemannian}
Let  $H : T^{*}M \longrightarrow \mathbb{R}$ be a Tonelli Hamiltonian and $L: TM\longrightarrow \bR$ the associated Tonelli Lagragian.
Let $W$ be a $C^{2}$ Lagrangian graph that is invariant by the Hamiltonian flow of $H$ and contained in a energy level $E(W)>c_{0}(L)$, where $c_{0}(L)$ denotes the strict Ma\~n\'e's critical value associated to $L$. Then, there exists a $C^{2}$ Riemannian metric $g$ on $M$, such that $W$ is an  invariant Lagrangian graph for the geodesic flow of $(M,g)$. 
\end{proposition}

\begin{proof}

The fact that $E(W)>c_{0}(L)$ ({\it i.e.}, it is supercritical) implies that the energy level containing $W$ admits a Finsler metric such that the Hamiltonian flow in the level is, up to reparametrization, the geodesic flow of the metric (Proposition \ref{propfinslerHam}). We can suppose without loss of generality that the energy level is $1$, and let $X_{F}$ be the geodesic vector field of the Finsler metric. The tautological 1-form $\alpha$ in the level is a contact form, and since $W$ is Lagrangian the restriction of the form $\alpha$ to $W$ is also closed. 

The geodesic vector field $X_{F}$ is the Reeb flow of the form $\alpha$, namely, it is the unique vector field such that $\alpha(X_{F}) = 1$ and $i_{X_{F}}^{*}d\alpha = 0$.

Since $W$ is  a Lagrangian graph, the canonical projection $\pi: T^*M \longrightarrow M$, $\pi(x,p):=x$, restricted to $W$ is a diffeomorphism and the push forward $\pi_{*}\alpha$ of the form $\alpha$ defines a closed 1-form on $M$. By the ideas of Subsection 2.5, we have a foliation $\mathcal{F}$ that is invariant by the flow $\psi_{t}: M \longrightarrow M$, obtained by projecting the Hamiltonian flow by the canonical projection. The foliation is defined by the closed 1-form $\pi_{*}\alpha$, in the sense that every leaf is always tangent to the kernel of the 1-form.  
Therefore, by choosing a Riemannian metric $(M,\bar{g})$ in $M$ we can apply Lemma \ref{Geodesible} to the 1-form $\pi^{*}\alpha$ and the vector field $X = d\pi (X_{F})$ tangent to the flow $\psi_{t}$; This concludes the proof of the proposition.
\end{proof}

\bigskip

Since the stable norms of any given two Finsler metrics on $M$ are equivalent, to get a cone in $H_{1}(M, \mathbb{R})$ it is enough to get a cone in the stable norm induced by the Riemannian metric $g$ in Proposition \ref{Riemannian}. This will be the goal of the forthcoming sections.

\subsection{The cone in homology containing the asymptotic homology of quasi-orbits}\label{secconehom}

Let $W \subset T^{1}\mathbb{T}^{n}$ be a $C^{2}$ Lagrangian graph invariant by the geodesic flow of $(\mathbb{T}^{n}, g)$. Let $\mathcal{F}$ be the foliation defined by the closed 1-form $\sigma$ associated to $W$, and let $\psi_{t}: \mathbb{T}^{n} \longrightarrow \mathbb{T}^{n}$ be the corresponding flow by geodesics orthogonal to the foliation. We know by Theorem \ref{tischler} that the closed 1-form defines a fibration of $\mathbb{T}^{n}$ over $\bS^{1}$. The main result of the section is the following:

\begin{theorem} \label{cone-fib}
Let $h \in H_{1}(\mathbb{T}^{n}, \mathbb{Z})$ be the homology class of the shortest closed geodesic of $(\mathbb{T}^{n}, g)$ that is not in the fundamental group of the fibers. Then, the asymptotic homologies generated by closed quasi-orbits of the flow $\psi_{t}$ are contained in the cone $\fC_{3.5}(h)$ with respect to the stable norm associated to $(\mathbb{T}^{n}, g)$. 
\end{theorem}

Let us recall the following definition.

\begin{definition}\label{quasiorbit}
A closed {\it quasi-orbit} with base point $x \in \mathbb{T}^{n}$ is a closed curve of the form 
$$c_{x,T} = (\cup_{t \in [0,T]}\psi_{t}(x)) \cup \Gamma_{x,T}$$
where $\Gamma_{x,T}$ is a minimizing geodesic of the metric $g$ joining $x$ and 
$\psi_{T}(x)$. The asymptotic homology generated by the quasi-orbits $c_{x,T}$ is the real homology class given by the limit (whenever it exists) of the homology classes 
$$ \lim_{T \rightarrow +\infty} \frac{1}{T}[c_{x,T}] $$
according to Definition \ref{defhomvector}. 
\end{definition}

{In the light of what we have discussed in Sections \ref{secfinslertoRiem} and \ref{sechamtoriem}, Theorem \ref{main} follows from the proof of Theorem \ref{cone-fib}.}\\

We shall divide the proof of Theorem \ref{cone-fib} into two steps.\\
Let $\tilde{\mathcal{F}}$, $\tilde{\psi}_{t}$ be respectively, the lifts of $\mathcal{F}$ and $\psi_{t}$ to the universal covering. We can suppose without loss of generality that the integral orbits of the flow $\psi_{t}$ are geodesics parametrized by arc length, so by Lemma \ref{Geodesible} the foliation $\mathcal{F}$ is preserved by the flow $\psi_{t}$. An analogous statement holds for $\tilde{\psi}_{t}$ and $\tilde{\mathcal{F}}$.\\

\noindent {\bf Case I - compact leaves.}

\noindent Let us first suppose that the leaves of $\mathcal{F}$ are all compact. By Corollary \ref{fibred-torus} we know that all the leaves of $\mathcal{F}$ defined a fibration of $\mathbb{T}^{n}$ over $\bS^{1}$ and that all the fibers are diffeomorphic to $\mathbb{T}^{n-1}$. 
\begin{proposition} \label{compactcase}
The asymptotic homologies generated by the quasi-orbits of $\psi_{t}$ are contained in the cone $\fC_{3}(h)$ in $H_{1}(\mathbb{T}^{n}, \mathbb{R})$ endowed with the stable norm associated to $(\mathbb{T}^{n}, g)$. 
\end{proposition}

\begin{proof}
Let $h \in H_{1}(\bS^{1},\mathbb{Z}) \subset H_{1}(\mathbb{T}^{n}, \mathbb{Z})$ be as in the statement of Theorem \ref{cone-fib}, and let $\gamma \subset \mathbb{T}^{n}$ be a closed geodesic in $h$ that minimizes the  length in the class $h$. Let $T_{\gamma} : \mathbb{R}^{n} \longrightarrow \mathbb{R}^{n}$ be the isometry of $(\mathbb{R}^{n}, \tilde{g})$ whose axes are the lifts of $\gamma$ to $\mathbb{R}^{n}$ by the covering map. 

Our main goal is to get an estimate for a cone in $H_{1}(\mathbb{T}^{n}, \mathbb{R})$ containing the asymptotic homologies generated by closed quasi-orbits of $\psi_{t}$. 

To do this, let us consider the lifted foliation $\tilde{\mathcal{F}}$ in the unversal covering $(\mathbb{R}^{n}, \tilde{g})$, and the lifted flow $\tilde{\psi}_{t}$. Let us choose $x \in \mathbb{T}^{n}$, $\tilde{x} \in \mathbb{R}^{n}$ such that $\Pi(\tilde{x}) = x$, where $\Pi : \bR^n \longrightarrow \mathbb{T}^{n}$ is the covering map. Let $\mathcal{F}(x)$ be the leaf of $\mathcal{F}$ containing $x$ and $\tilde{\mathcal{F}}(\tilde{x})$ the corresponding lifted leave. Let us consider the action of $T_{\gamma}$ on $\tilde{\mathcal{F}}(\tilde{x})$. Notice that the leaves are preserved by the action of $T_{\gamma}$, {\it i.e.}, 
$$ T_{\gamma}(\tilde{\mathcal{F}}(\tilde{x})) = \tilde{\mathcal{F}}(T_{\gamma}(\tilde{x})) $$
for every $\tilde{x} \in \tilde{\mathbb{R}^{n}}$. The leaves $\tilde{\mathcal{F}}$ are equidistant, so let $\lambda $ be the distance from $\tilde{F}(\tilde{x})$ to $T_{\gamma}(\tilde{F}(\tilde{x}))= \tilde{\psi}_{\lambda}(\tilde{\mathcal{F}}(\tilde{x})$. 
\medskip

We clearly have a return map $P_{0}:\mathcal{F}(x) \longrightarrow \mathcal{F}(x)$ for the leaf $\mathcal{F}(x)$ associated to $T_{\gamma}$. Namely, $P_{0}$ is given by the following commuting relation:
$$ \Pi (\tilde{\psi}_{m\lambda}(\tilde{p})) = P_{0}^{m}(\Pi(\tilde{p}))$$ 
for every $\tilde{p} \in \tilde{\mathcal{F}}(\tilde{x})$. 

Let $d$ be the distance in $(\mathbb{R}^{n}, \tilde{g})$ (we shall use by abuse of notation the same notation for the distances in $(\mathbb{T}^{n}, g)$ and $(\mathbb{R}^{n}, \tilde{g})$, the metric space will determine uniquely the corresponding distance).  
\medskip

\noindent \textbf{Claim 1}: {\it Let $\ell_{g}(\gamma)$ be the length of $\gamma$, and let $d_{m}=m\lambda$ be the distance from $\tilde{\mathcal{F}}(\tilde{x})$ to $T_{\gamma}^{m}(\tilde{\mathcal{F}}(\tilde{x}))$. Then $d_{m} \leq ml_{g}(\gamma)$ for every $m \in \mathbb{N}$.} 
\medskip

Indeed, suppose without loss of generality that a lift $\tilde{\gamma}$ of $\gamma$ is such that $\tilde{\gamma}(0) \in \tilde{F}(\tilde{x})$. Since $T_{\gamma}$ preserves $\tilde{\gamma}$ and  
$$d(\tilde{\gamma}(t), T_{\gamma}^{m}(\tilde{\gamma}(t)) ) = d(\tilde{\gamma}(t) , \tilde{\gamma}(t +m\ell_{g}(\gamma)) ) = m\ell_{g}(\gamma) $$
we have that $m\ell_{g}(\gamma)$ is the length of a curve joining two points in $\tilde{\mathcal{F}}(\tilde{x})$ and $T_{\gamma}^{m}(\tilde{\mathcal{F}}(\tilde{x}))$, that is another lift of $\mathcal{F}(x)$. Since the leaves of $\tilde{\mathcal{F}}$ are equidistant, and $T_{\gamma}$ acts by isometries, we have that the distance from $\tilde{\mathcal{F}}(\tilde{x})$ to  
$T_{\gamma}^{m}(\tilde{\mathcal{F}}(\tilde{x}))$ is $m\lambda$. Moreover, $m\lambda$ coincides with the length of the minimizing geodesic $\tilde{\psi}_{t}(\tilde{p})$, for $\tilde{p} \in \tilde{\mathcal{F}}(\tilde{x})$ and $t \in [0,m\lambda]$. Thus, the distance $m\lambda$ must be less than or equal to the length of any curve joining a point in $\tilde{\mathcal{F}}(\tilde{x})$ to a point in  
$T_{\gamma}^{m}(\tilde{\mathcal{F}}(\tilde{x})$, so $d_{m}= m\lambda \leq ml_{g}(\gamma)$ as claimed. 
\medskip

Let $y_{\tilde{x}} \in \tilde{\mathcal{F}}(\tilde{x}) $ be the closest point to $\tilde{x}$ in an axis of $T_{\gamma}$, and let $\hat{\gamma}$ be this axis, with the property that $\hat{\gamma}(0) = y_{\tilde{x}}$. 
Let us adopt the notation $[[z,w]]$ for a minimizing geodesic in $(\mathbb{R}^{n}, \tilde{g})$ joining $z, w$.

Consider the following quasi-orbits in $\mathbb{T}^{n}$:
$$ c_{x,m} = \psi_{[0,m\lambda]}(x)\cup [[\psi_{m\lambda}(x),x]]$$
for $ m \in \mathbb{N}$, where $\psi_{[0,m\lambda]}(x) = \cup_{t \in [0,m\lambda]}\psi_{t}(x)$. 

Let us consider the lifts of the curves $c_{x,m}$ with base point $\tilde{x}$: 
$$ \tilde{c}_{\tilde{x},m} = \tilde{\psi}_{[0,m\lambda]}(\tilde{x}) \cup [[\tilde{\psi}_{m\lambda}(\tilde{x}), \tilde{x}_{m}]]$$ 
where the geodesic $[[\tilde{\psi}_{m\lambda}(\tilde{x}), \tilde{x}_{m}]]$ is the lift of the geodesic $[[\psi_{m\lambda}(x),x]]$ starting at $\tilde{\psi}_{m\lambda}(\tilde{x})$, and $\Pi(\tilde{x}_{m}) = x$. 

Let us consider the curves in $\mathbb{R}^{n}$
\begin{eqnarray*}
 \tilde{\sigma}^{1}_{\tilde{x},m} & = & [[\tilde{x}, y_{\tilde{x}}]] \cup \hat{\gamma}[0,m\ell_{g}(\gamma)]
\cup[[\hat{\gamma}(m\ell_{g}(\gamma)), T_{\gamma}^{m}(\tilde{x})]] \\
 & = & [[\tilde{x}, y_{\tilde{x}}]] \cup [[ y_{\tilde{x}}, T_{\gamma}^{m}(y_{\tilde{x}}) ]] \cup [[T_{\gamma}^{m}(y_{\tilde{x}}),  T_{\gamma}^{m}(\tilde{x})]]
\end{eqnarray*}
and 
$$\tilde{\sigma}^{2}_{\tilde{x},m} =  [[ T_{\gamma}^{m}(\tilde{x}), \tilde{x}_{m}]].$$

Let $$\tilde{\Gamma}_{\tilde{x},m} = \tilde{\sigma}^{1}_{\tilde{x},m} \circ \tilde{\sigma}^{2}_{\tilde{x},m}$$
where $a\circ b$ is the concatenation of the curves $a$, $b$. 

The images by the covering map 
$$\sigma^{1}_{x,m}= \Pi(\tilde{\sigma}^{1}_{\tilde{x},m}), \mbox{ }\sigma^{2}_{x,m}= \Pi(\tilde{\sigma}^{2}_{\tilde{x},m}), \mbox{ } 
\Gamma_{x,m} = \Pi(\tilde{\Gamma}_{\tilde{x},m})$$ 
are closed curves with base point $x$. Moreover, the curve 
$$ \tilde{c}_{\tilde{x},m} \circ \tilde{\Gamma}_{\tilde{x},m} $$
is a closed curve in $\mathbb{R}^{n}$. Therefore, the curves $c_{x,m}$ and $\Gamma_{x,m}$ have the same homology class and we have 

$$ [c_{x,m}] =  [\Gamma_{x,m}] =   [\sigma^{1}_{x,m}] + [\sigma^{2}_{x,m}].$$ 
Since the homology class of $\sigma^{1}_{x,m}$ is the homology class of $\gamma[0,m\ell_{g}(\gamma)]$, we have 
$$ [c_{x,m}] =  [\Gamma_{x,m}] =   [h^{m}]+ [\sigma^{2}_{x,m}].$$ 
With this decomposition of the homology class $[ \Gamma_{x,m} ]$, we can estimate its stable norm in 
terms of the stable norm of $h= [ \gamma ]$. 
\bigskip

\begin{figure}[ht]
\begin{center}
\includegraphics[scale=0.4]{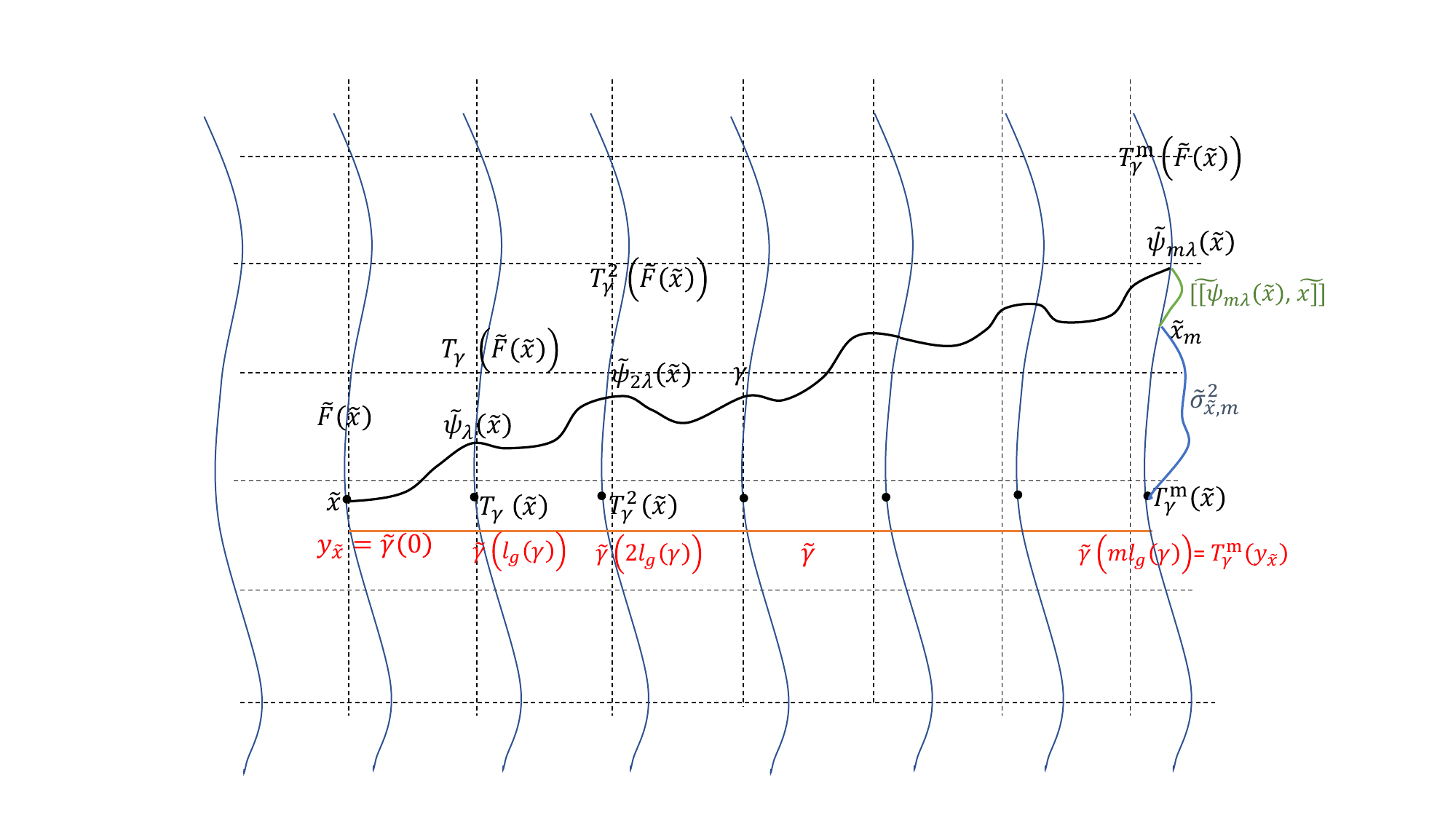}
\begin{picture}(0,0)
\end{picture}
\end{center}
\vspace{-0.5cm}
\caption{{\small The lift of a quasi-orbit with base point $\tilde{x}$.}}\label{fig-1}
\end{figure}

\bigskip

First of all, the stable norm of $\sigma^{1}_{x,m}$ with respect to the metric $g$ is 
$$ \| \sigma^{1}_{x,m} \|_{s} = \| h^{m} \|_{s} = m \| h \|_{s} = m\ell_{g}(\gamma)  .$$

To estimate the stable norm of $[\sigma^{2}_{x,m}]$, observe that the length of $\sigma^{2}_{x,m}$ coincides with the length of $\tilde{\sigma}^{2}_{\tilde{x},m}$, that is a minimizing geodesic joining $\tilde{x}_{m}$ and $T_{\gamma}^{m}(\tilde{x})$.  By the triangle inequality, we have
$$ d(\tilde{x}_{m}, T_{\gamma}^{m}(\tilde{x})) \leq\ell_{\tilde{g}}(\tilde{c}_{\tilde{x},m}) +\ell_{\tilde{g}}(\tilde{\sigma}^{1}_{\tilde{x},m}) $$
By Claim 1, the length of $\tilde{c}_{\tilde{x},m}$ satisfies
	$$\ell_{g}(\tilde{c}_{\tilde{x},m}) \leq ml_{g}(\gamma) + D  $$ 
where $D$ is the diameter of $(\mathbb{T}^{n}, g)$. 
Moreover, since the length of $\tilde{\sigma}^{1}_{\tilde{x},m}$ satisfies 
$$\ell_{g}(\tilde{\sigma}^{1}_{\tilde{x},m} ) \leq 2D + m\ell_{g}(\gamma) $$ 
we conclude that 
$$\ell_{g}([\tilde{x}_{m}, T_{\gamma}^{m}(\tilde{x})]) \leq 2m\ell_{g}(\gamma) + 3D .$$
Therefore we get that 
$$ \| \sigma^{2}_{x,m} \|_{s} \leq 2 \| h^{m} \|_{s} + 3D  .$$
Finally, we deduce that the stable norm in homology of the closed curve $c_{x,m}$ satisfies
\begin{eqnarray*}
 \| c_{x,m}) \|_{s} =  \| \Gamma_{x,m}  \|_{s} & \leq & \| h^{m} \|_{s} + 2 \| h^{m} \|_{s} + 3D  \\
& = & 3 \| h^{m} \|_{s} + 3D  
\end{eqnarray*}
that implies 
$$ \frac{\| c_{x,m}) \|_{s}}{ \| h^{m} \|_{s}} \leq 3 + \frac{3D}{m\ell_{g}(\gamma)}. $$
and hence
$$\limsup_{m \rightarrow \infty} \frac{\| c_{x,m}) \|_{s}}{ \| h^{m} \|_{s}} \leq 3 . $$
This clearly yields that every convergent subsequence of homologies of the quasi-orbits $c_{x,m}$ must be contained in the cone $\fC_{3}(h)$, as claimed. \end{proof}
\bigskip

\noindent {\bf Case II - Noncompact leaves case: geodesible flows.} 
\noindent

\begin{lemma} \label{noncompact}
Suppose that the leaves of $\mathcal{F}$ are not compact. Then the asymptotic homologies generated by the quasi-orbits of $\psi_{t}$ are contained in the cone 
$\fC_{3.5}(h)$ in $H_{1}(\mathbb{T}^{n}, g)$ endowed with the stable norm of $(\mathbb{T}^{n}, g)$. 
\end{lemma}

\begin{proof}
If the leaves are non compact, given $\epsilon >0$ we can apply Corollary \ref{fibred-torus} to $\epsilon$-approximate the closed 1-form $\alpha$ by a closed 1-form $\alpha_{\epsilon}$ in the $C^{2}$ topology, defining a fibration $\mathcal{F}_{\epsilon}$ by tori in $\mathbb{T}^{n}$. The flow $\psi_{t}$ is ``almost perpendicular'' to the fibers of $\mathcal{F}_{\epsilon}$, and then, Lemma \ref{Geodesible} implies that 
\bigskip

\textbf{Claim :} 
{\it There exist $\delta(\epsilon)$ depending continuously on $\epsilon $, with $\delta(0) = 0$, and a Riemannian metric $g_{\epsilon}$ in $\mathbb{T}^{n}$ that is $\delta(\epsilon)$-$C^{2}$ close to $g$, with the following properties:
 
\begin{enumerate} 
	\item The orbits of $\psi_{t}$ are geodesics of $(\mathbb{T}^{n}, g_{\epsilon})$. 
	\item The fibers of $\mathcal{F}_{\epsilon}$ are perpendicular to the orbits of $\psi_{t}$. 
\end{enumerate}}

The fibers of $\mathcal{F}_{\epsilon}$ are perpendicular to the vector field $X$ with respect to $g_{\epsilon}$, and the orbits of $\psi_{t}$ are geodesics of $g_{\epsilon}$ by Lemma \ref{Geodesible}. This shows the Claim. 
\medskip

Thus, Proposition \ref{compactcase}
implies that the asymptotic homologies generated by $\psi_{t}$ are contained in the cone $\fC_{3}^{\epsilon}(h)$ in $H_{1}(\mathbb{T}^{n}, g)$ endowed with the stable norm of $(\mathbb{T}^{n}, g_{\epsilon})$. 

Since the stable norm depends continuously on the Riemannian metric, the cones 
$\fC_{3}^{\epsilon}(h)$ will approach the cone $\fC_{3}(h)$ uniformly on compact sets. In particular, we can choose $\epsilon$ small enough in a way that 
$$\fC_{3}^{\epsilon}(h) \subset \fC_{3.5}(h).$$ 
This finishes the proof of the Lemma. 
\end{proof}
\medskip

The proof of Theorem \ref{cone-fib} follows from Proposition \ref{compactcase} and Lemma \ref{noncompact}. This concludes the proof of Theorem \ref{main}.\\

\begin{remark}\label{remmechanical} {\bf (Mechanical systems and geodesic flows)}
It follows from Proposition \ref{compactcase} that the size of the cone in Theorem \ref{main} can be chosen uniformly in the category of Riemannian metrics. The same holds for mechanical Hamiltonian systems, by applying the classical {\it Maupertuis principle}. In fact, for a compact Riemannian manifold, $(M, g)$, let
$$
    H(x, p) := \frac{1}{2} \sum_{ij} g^{ij}(x)p_ip_j - V(x)
$$
be a  mechanical Hamiltonian function on $T^*M$, where $V \in C^2(M)$ denotes some potential function. This is a Tonelli Hamiltonian and one can check that Ma\~n\'e's strict critical value equals $- \min_x V(x)$ (see for example \cite{Gonzalonotes, Sorrentinobook}).
For every $e > - \min_x V(x)$, 
the energy level $\Sigma_e = \{H(x, p) = h\}$ is also an energy level for the Hamiltonian system given by 
$$
    H_e(x, p) := \frac{1}{2} \sum_{ij} g^{ij}(x) \left( h + V(x) \right) p_ip_j.
$$
Maupertuis principle (see for instance \cite{BKF}) states that the integral curves for the Hamiltonian vector fields $X_H$ and $X_{H_e}$ on the fixed energy level $\Sigma_e$ coincide. 
Note that the vector field $X_{H_e}$ gives rise to the geodesic flow of the Riemannian metric $g_e$ with
$$
    g_e^{ij}(x) = (h + V(x)) g^{ij}(x).
$$
\end{remark}

\bigskip

\subsection{A counterexample to Theorem \ref{main} at the critical value}\label{maneexample}
In this section we describe  an example to show that Theorem \ref{main} might fail on the energy level corresponding to Ma\~n\'e's strict critical value.

Let $M=\T^2$ equipped with the flat metric and consider a vector field $X$ with norm $1$ and such that X has two closed orbits $\g_1$ and $\g_2$ in opposite (non-trivial) homology classes and any other orbit asymptotically approaches $\g_1$ in forward time and $\g_2$ in backward time; for example one can consider $X(x_1,x_2)=(\cos(2\pi x_1), \sin(2\pi x_1))$, where $(x_1,x_2)\in \T^2 = \R^2/\Z^2$.

\begin{figure}[h!]
\begin{center}
\includegraphics[scale=0.2]{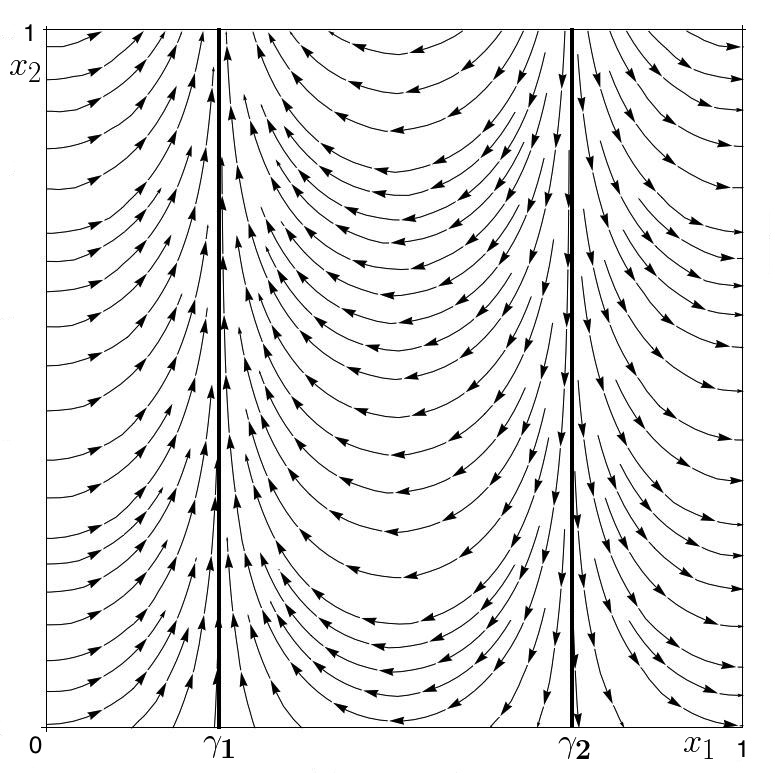}
\caption{Plot of the vector field $X$.}
\end{center}
\end{figure}

One can embed this vector field into the Euler-Lagrange vector field given by the Tonelli Lagrangian 
$L_X(x,v)= \frac{1}{2}\|v-X(x)\|^2$ (these are sometimes called {\it Ma\~n\'e's Lagrangians}, see  \cite{Mane1, Gonzalonotes, Sorrentinobook}).  It is easy to check that the corresponding Hamiltonian is given by $H_X(x,p)= \frac{1}{2}\|p\|^2 + \langle p, X(x) \rangle$. In particular, the zero section $\Lambda_0 := \T^2\times \{0\}$ is invariant under the flow of $H_X$ and the flow restricted to $\Lambda_0$ is conjugate to the flow of the vector field $X$.
Hence, $\Lambda_0$ contains two periodic orbits of opposite asymptotic homologies and all other orbits are heteroclinic to these periodic orbits. This implies that the set of asymptotic homologies cannot be contained in any proper cone. See also \cite[Remark 3.3.5, iii)]{Sorrentinobook} for more details.\\

\medskip

\section{Rotation vector of diffeomorphisms of the 2-torus and asymptotic homology of 3-dimensional Lagrangian invariant graphs} \label{secdim3}

In this section we want now to describe how to apply  Theorem \ref{main} to give a partial answer to the following question: \\

\noindent {\sc Question:}  {\it If the smallest cone in homology generated by  the asymptotic homologies of quasi-orbits  of an invariant Lagrangian graph has non-empty interior, then has the flow a periodic orbit?} \\

An affirmative answer to this question would imply that if an invariant Lagrangian graph contains no periodic orbits then the asymptotic homologies of quasi-orbits of the Lagrangian graph are all contained in a plane of nontrivial codimension. This is interesting from the point of view of Aubry-Mather theory and  differentiability, or lack thereof, of {\it Mather's minimal average action functions}, the so-called $\alpha$, $\beta$ functions (see \eqref{defalpha} and \eqref{defbeta} for a definition of these functions). \\

We shall now prove Theorem \ref{main2}: its proof is based on the relationship between the rotation set of a homeomorphism of $\mathbb{T}^{2}$ that is homotopic to the identity, and the asymptotic homologies of quasi-orbits of its suspension flow in $\mathbb{T}^{3}$. The rotation set of a homeomorphism of the two torus homotopic to the identity was defined and studied by J. Franks  (see \cite{kn:Franks}). 

\subsection{Invariant Lagrangian graphs and suspensions} \label{sec5-1}

Let us start exploring some of the consequences of the results of Section 4. We begin by recalling the definition of a {\it suspension flow}. 

\begin{definition}
Let $M$ be a metric space and let $f : M \longrightarrow M$ be a homeomorphism. Let us define an equivalence relation $\sim$ between $M\times [0,1]$ and $M \times [0,1]$ by $(x,1) \sim (f(x),0)$, and the equivalence class of $(x,t)$ is formed by a single element if $t \neq \{0, 1\}$. The suspension flow of $f$ is the flow $f_{t} : M\times [0,1]/\sim \longrightarrow M\times[0,1]/\sim $ given by 
$$ f_{t}(x,s) = \psi_{t+s}(x,0) = (f^{n}(x), t+s \; ({\rm mod}\; 1)) $$
where $n$ is the integer part of $t+s$. \\
\end{definition}

\smallskip

\begin{lemma} \label{suspension}
Let $H: T^{*}\mathbb{T}^n \longrightarrow \bR$ a Tonelli Hamiltonian on the torus and let $W \subset T^{*}\mathbb{T}^{n}$ be a $C^{2}$ Lagrangian graph invariant by the Hamiltonian flow contained in a supercritical energy level $E > c_{0}(L)$, where $L$ is the corresponding Lagrangian. 
Let $\psi_{t} :  \mathbb{T}^n \longrightarrow \mathbb{T}^n$ be the canonical projection of the Hamiltonian flow restricted to $W$. Then $\psi_{t}$ is conjugate to the suspension of a homeomorphism $f : \mathbb{T}^{n-1} \longrightarrow \mathbb{T}^{n-1}$. 
\end{lemma}

\begin{proof}

By the results of Section 4 we know that there exists a Riemannian metric $g$ in $\mathbb{T}^{n}$ such that the orbits of the flow $\psi_{t}$ are geodesics of $(\mathbb{T}^{n}, g)$ which are perpendicular to a fibration by compact leaves, all diffeomorphic to $\mathbb{T}^{n-1}$. 

 Parametrizing the orbits of $\psi_{t}$ by arclength, we have that the parameter $t$ of the flow is locally the distance between the leaves of the fibration, moreover, the leaves are preserved by $\psi_{t}$. Let $h \in H_{1}(\mathbb{T}^{n}, \mathbb{Z})$ be the shortest element that is not in the first homology group of the leaves.  

As in the proof of Proposition \ref{compactcase}, given a leaf $\mathcal{F}_{0}$ of the fibration, we have a return map $P_{0}:\mathcal{F}_{0} \longrightarrow \mathcal{F}_{0}$ associated to $h$, namely, there exists $\lambda >0$ depending on $h$ such that 
$$ P_{0}(x) = \psi_{\lambda}(x)$$ 
for every $x \in \mathcal{F}_{0}$. Let $f_{t} :\mathbb{T}^{n} \longrightarrow \mathbb{T}^{n}$ be defined by $f_{t}(x) = \psi_{\lambda t }(x)$. It is clear that 
$$f_{m}(x)  = P_{0}^{m}(x) \in \mathcal{F}_{0}$$ 
for every $x \in \mathcal{F}_{0}$.  It is easy to check that $f_{t}$ is conjugate to the suspension of  $P_{0}$ and that it is conjugate to $\psi_{t}$. 
\end{proof}

\medskip

\subsubsection{Rotation set of a homeomorphism homotopic to the identity}

We follow \cite{kn:Franks}. Let $f : \mathbb{T}^{2} \longrightarrow \mathbb{T}^{2}$ be a continuous map homotopic to the identity. 
Let $F : \mathbb{R}^{2} \longrightarrow \mathbb{R}^{2}$ be a lift of $f$ to the universal covering of $\mathbb{T}^{2}$. The rotation set $\rho(F)$ is the set of accumulation points of the vectors 
$$ \left\{ \frac{F^{n}(x) - x}{n} , \mbox{ } x \in \mathbb{R}^{2}, \mbox{ } n \in \mathbb{N}^{+} \right\}. $$ 

It was shown by Misiurewicz--Ziemian \cite{kn:MZ} that when $f$ is a homeomorphism, the rotation set $\rho(F)$ is a convex set. 

The theory of rotation sets extends in many senses Poincar\'{e}'s theory of rotation vectors for homeomorphims of the circle. A remarkable 
result due to Franks \cite{kn:Franks} is the following:

\begin{theorem}\label{Franks}
	Let $f : \mathbb{T}^{2} \longrightarrow \mathbb{T}^{2}$ be a homeomorphism homotopic to the identity and let $F: \bR^2\longrightarrow \bR^2$ denote its lift. If a vector $v$ with rational coordinates lies in the interior of the rotation set $\rho(F)$, then the homeomorphism $f$ has a periodic point. 
\end{theorem}

The definitions of rotation set and asymptotic homology are quite similar. Indeed, let $\Pi : \mathbb{R}^{2} \longrightarrow \mathbb{T}^{2}$  be the covering map. Let us endow $\mathbb{T}^{2}$ with a Riemannian metric $g$, and consider a homeomorphism $f : \mathbb{T}^{2} \longrightarrow \mathbb{T}^{2}$ and $x \in \mathbb{R}^{2}$. Let us lift the metric $g$ to a metric $\tilde{g}$ in the universal covering of $\mathbb{T}^{2}$, and consider the curves $C_{n}(x)$ formed by the union of a minimizing geodesic $[x, F^{n}(x)]$ joining $x$ to $F^{n}(x)$, and a minimizing geodesic $[F^{n}(x), x_{n}]$ joining $F^{n}(x)$ to a lift $x_{n}$ of $x$ in the same fundamental domain of $F^{n}(x)$. The curve 
$$c_{n}(\Pi(x) ) = \Pi([x, F^{n}(x)] \cup [F^{n}(x), x_{n}])$$ 
is a closed continuous curve with base point $\Pi(x)$. Let $[c_{n}(\Pi(x))]$ be the integer homology class of $c_{n}(\Pi(x))$. Then we have 

\begin{lemma} \label{rotation-ashom}
	The rotation set $\rho(F)$ coincides with the set of accumulation points of 
	$$ \left\{ \frac{[c_{n}(\Pi(x))] }{n},\mbox{ } x \in \mathbb{R}^{2}, \mbox{ } n \in \mathbb{N}^{+} \right\}.$$
\end{lemma}

\begin{proof}

The homology class $[c_{n}(\Pi(x))]$ is represented by the integer vector 
$$x_{n} - x = (x_{n} - F^{n}(x) ) + (F^{n}(x) - x) .$$
Notice that the diameter $D$ of a fundamental domain of $\mathbb{T}^{2}$ is an upper bound for the lenght of the vectors $(x_{n} - F^{n}(x) )$. So we get 

$$ \frac{[c_{n}(\Pi(x))] }{n} = \frac{(x_{n} - F^{n}(x) )}{n}  + \frac{(F^{n}(x) - x) }{n} $$ 
from which it is clear that any accumulation point of the set $\left\{ \frac{[c_{n}(\pi(x))]}{n} \right\}$ is an accumulation point of the set $ \frac{(F^{n}(x) - x)}{n}$, and vice-versa.  This proves the lemma.

	\end{proof}

\subsubsection{Asymptotic homologies of a suspension flow}
 
Let us consider a $C^{2}$ diffeomorphism $f : \mathbb{T}^{2} \longrightarrow \mathbb{T}^{2}$ homotopic to the identity, and its suspension flow 
$\Psi_{t} : \mathbb{T}^{3} \longrightarrow \mathbb{T}^{3}$. 

The domain $\mathbb{T}^{3}$ of the suspension is naturally fibered over $\bS^{1}$ with fiber $\mathbb{T}^{2}$. The fibers, $\{ t \} \times \mathbb{T}^{2}$ for $t \in [0,1]$, are all transversal to the vector field tangent to $\Psi_{t}$. So to each $C^{2}$ family of metrics $g_{t}$ defined in $\{ t \} \times \mathbb{T}^{2}$ for $t \in [0,1]$, with $g_{0}= g_{1}$, we can associate a $C^{2}$ metric $g$ in $\mathbb{T}^{3}$ such that the orbits of the flow $\Psi_{t}$ are geodesics of this metric, by the same procedure described in Subsection 4.1. The homology group of $\mathbb{T}^{3}$ is the direct sum 
$$ H_{1}(\mathbb{T}^{3}, \mathbb{R}) =  \langle h\rangle \oplus H_{1}(\mathbb{T}^{2}, \mathbb{R}) $$
where $h $ is the generator of the basis of the fibration in integer homology. Lemma \ref{rotation-ashom} and the topological arguments in the proof of Theorem \ref{cone-fib} imply that 

\begin{lemma} \label{product-hom}
	The asymptotic homology generated by quasi-orbits of the flow $\Psi_{t}$ is the set of vectors
$$\{ (h, \sigma):\; \sigma \in \rho(F) \}$$
where $F$ is a lift of $f$ in $\mathbb{R}^{2}$.
\end{lemma}

\begin{proof}
	
Let $h \in H_{1}(\mathbb{T}^{3}, \bR)$   be the homology class of the basis of the fibration induced by the suspension of $f$. Let $\tilde{x}\in \mathbb{R}^{3}$ be a base point in the universal covering of $\mathbb{T}^{3}$, and let $\Pi: \mathbb{R}^{3} \longrightarrow \mathbb{T}^{3}$ be the covering map. Let $\tilde{g}= \Pi^{*}(g)$ be the pullback of the Riemannian metric $g$ to the universal covering. 

Let $T : \mathbb{R}^{3} \longrightarrow \mathbb{R}^{3}$ be the covering $\tilde{g}$-isometry associated to the element $h$ viewed as an element of $\pi_{1}(\mathbb{T}^{3})$. The fibration $\mathcal{F}$ of $\mathbb{T}^{2}$ lifts to a foliation $\tilde{\mathcal{F}}$ of the universal covering by leaves  $\tilde{\mathcal{F}}(x)$, meaning the leaf containing a point $x \in \mathbb{R}^{3}$. The foliation is preserved by the isometry $T$, namely, $T(\tilde{\mathcal{F}}(x) )= \tilde{\mathcal{F}}(T(x))$ for every $x \in \mathbb{R}^{3}$. 

Let $F: \tilde{\mathcal{F}}(\tilde{x}) \longrightarrow \tilde{\mathcal{F}}(\tilde{x})$ be a lift of the diffeomorphism $f$. Let $\tilde{\phi}_{t} : \mathbb{R}^{3} \longrightarrow \mathbb{R}^{3}$ be a lift of the suspension flow $\phi_{t}$. 

The definition of the suspension implies at once that 
$$ T(F(x)) = \tilde{\phi}_{1}(x)$$ 
for every $x \in \tilde{\mathcal{F}}(\tilde{x})$. In particular, $T^{n}(F^{n}(x)) = \tilde{\phi}_{n}(x)$ for every $x \in \tilde{\mathcal{F}}(\tilde{x})$ and 
$n \in \mathbb{Z}$. 

	Now, let $x \in \mathbb{T}^{3}$. 	The asymptotic homology generated by the quasi-orbits with base point $x$ is the set of accumulation points of the homologies of the curves 
	$$ C_{n} = \phi_{[0,t_{n}]}(x) \cup [[\phi_{t_{n}}(x) , x]]$$ 
	where $[[\phi_{t_{n}}(x) , x]]$ is a minimizing geodesic joining $\phi_{t_{n}}(x)$ to $x$. 
	
	Let $t_{n} \in [m_{n}, m_{n}+1]$, where $m_{k}$ is an integer for every $k \in \mathbb{Z}$. It is easy to show that the set of accumulation points of these curves is the same as the set of accumulation points of the curves 
	$$ \mathcal{C}_{n} = \phi_{[0,m_{n}]}(x) \cup[[\phi_{m_{n}}(x) , x_{n}]]_{\mathcal{F}}$$ 
	where $x_{n} \in \mathcal{F}(x)$ is close to $x$, and $[[\phi_{m_{n}}(x) , x_{n}]]_{\mathcal{F}}$ is a minimizing geodesic joining $\phi_{m_{n}}(x)$ to $x_{n}$ in the fiber $F(x)$,  endowed with the restriction of the metric $g$. In other words, we can restrict the analysis of the asymptotic homology of $\phi_{t}$ to the analysis of the asymptotic homology of the flow orbits of the Poincar\'{e} map of the fiber $\mathcal{F}(x)$ induced by the flow, namely, the diffeomorphism $f$. 
	
	Now, we can argue as in the proof of Theorem \ref{cone-fib}. The curve $\mathcal{C}_{n}$ lifts to a curve 
	$$\tilde{\mathcal{C}}_{n} = \tilde{\phi}_{m_{n}}(\tilde{x}) \cup [[\tilde{\phi}_{m_{n}}(\tilde{x}), \tilde{x}_{n}]]_{\tilde{\mathcal{F}}}$$ 
    in the universal covering, joining $\tilde{x} \in \tilde{\mathcal{F}}(\tilde{x})$ to $\tilde{x}_{n} \in \tilde{\mathcal{F}}(\tilde{\phi}_{m_{n}}(\tilde{x}))= T^{m_{n}}(\tilde{\mathcal{F}}(\tilde{x}))$. 
    The curve can be written as the composition 
    $$ \tilde{\mathcal{C}}_{n} = \Gamma_{n} \cup [[T^{m_{n}}(\tilde{x}), \tilde{x}_{m_{n}}]]_{\tilde{\mathcal{F}}}$$
  
    where $\Gamma_{n}$ is formed by the union of minimizing geodesics $\gamma_{i}= T^{i}(\gamma)$, $0 \leq i \leq m_{n}$, $\gamma$ is a minimizing geodesic joining $\tilde{x}$ to $T(\tilde{x})$, and $[[T^{m_{n}}(\tilde{x}), \tilde{x}_{m_{n}}]]_{\tilde{\mathcal{F}}}$ a minimizing geodesic  in the leaf $\tilde{\mathcal{F}}(\tilde{x}_{n})$, with respect to the restriction of the metric $\tilde{g}$ to the leaf.

    $$ [  \Pi(\mathcal{C}_{n})] = ( m_{n} h , \rho_{n})$$ 
    Finally, the accumulation points of the classes $\frac{1}{m_{n}}[  \Pi(\mathcal{C}_{n})] $ are of the form $( h , \sigma)$, where $\sigma $ is in the rotation set of $f$ according to Lemma \ref{rotation-ashom}, as we wished to show.

\end{proof}

Since $\rho(F)$ is compact, there exists $L>0$ such that 
$$ \frac{ \parallel \sigma \parallel_{s} }{\parallel h \parallel_{s} } \leq L $$
for every $\sigma \in \rho(F)$. This implies that the set of vectors $(th , t \sigma)$, $t \in \mathbb{R}^{+}$, $\sigma \in \rho(F)$, generates a cone $\fC_{L}(h)$ in $H_{1}(\mathbb{T}^{3}, \mathbb{R})$. So we get an alternative proof of Theorem \ref{cone-fib} for this particular case of suspension flows of diffeomorphims in $\bT^{2}$ that are homotopic to the identity. 

Moreover, the rotation set $\rho(F)$ is a convex set according to \cite{kn:MZ} contained in $H_{1}(\mathbb{T}^{2}, \mathbb{R})$, that is identified with $\mathbb{R}^{2}$. Thus, the cone $\fC_{L}(h)$ is a convex set in $H_{1}(\mathbb{T}^{3}, \mathbb{R})$. 
\bigskip

\subsection{Proof of Theorem \ref{main2}} \label{sec5-2}

We can now prove the main result of this section.

\begin{proof} 
Under the assumptions of Theorem \ref{main2} suppose by contradiction that the cone $\fC_{L}(h)$ generated by the asymptotic homologies of the suspension flow $\Psi_{t}$ is not contained in any codimension-$1$ plane in $H_{1}(\mathbb{T}^{3}, \mathbb{R})$. Then there exist 3 linearly independent vectors in $\fC_{L}(h)$.

 By Lemma \ref{product-hom}, there exists 3 linearly independent vectors 
$$ (h, \sigma_{1}), \mbox{ } (h, \sigma_{2}), \mbox{ }(h, \sigma_{3}) \in H_{1}(\mathbb{T}^{3}, \mathbb{R})= \langle h \rangle \oplus H_{1}(\mathbb{T}^{2}, \mathbb{R}),$$ 
where each $\sigma_{i}$ is an element of the rotation set of the homeomorphism $f$. This implies that the points $\sigma_{i}$ are in general position in the rotation set of $f$. Since the rotation set of $f$ is a convex set of the plane $H_{1}(\mathbb{T}^{2}, \mathbb{R})$, this yields that its interior is non-empty. By Franks' Theorem, $f$ has periodic points. Since each periodic orbit of the flow $\Psi_{t}$ is the suspension of a periodic orbit of $f$, we conclude that the flow $\Psi_{t}$ has a periodic orbit, contradicting the hypothesis of the Theorem. 
\end{proof}

\medskip

\begin{remark}
The main obstructions to extend the proofs of Theorems \ref{cone-fib} and \ref{main2} to any Lagrangian graph are mainly two: 
\begin{enumerate}
	\item If the dimension of the Lagrangian graph is $3$, the Poincar\'{e} return map might not be homotopic to the identity in the cross section of dimension 2, so the notion of rotation set is not well defined. 
	\item In higher dimensions the rotation set is not well defined for homeomorphisms which are not homotopic to the identity, and Franks' Theorem \ref{Franks} is not available. 
\end{enumerate}
\end{remark}

\medskip

\subsection{An algebraic criterion for the differentiability of Mather's $\alpha$ function}\label{sec5-3}

We would like to explore under which assumptions in the asymptotic homology of an invariant Lagrangian graph the so-called Mather's $\alpha$ function (see \eqref{defalpha}) is differentiable. The results of the previous subsection show that the absence of periodic orbits of the projection of the Hamiltonian flow restricted to a Lagrangian graph implies that the asymptotic homologies generated by quasi-orbits are contained in a codimension one plane in the first homology group. This is somehow related to the differentiability of the Mather's $\alpha$ function,  by the following reasons:
\begin{enumerate}
\item The quasi-orbits associated to generic (from the measure point of view) orbits in the support of a Mather measure have a (unique) naturally associated asymptotic homology, given by the Schartzmann cycle. 
\item Every invariant Lagrangian graph has Mather measures supported in the graph. 
\item If all Mather measures supported in the Lagrangian graph have the same asymptotic homology, then the Mather's $\alpha$ function is differentiable (see the Appendix). 
\end{enumerate}

Therefore, if the minimal cone in homology containing the asymptotic homologies of quasi-orbits in the support of Mather measures in the graph reduces to a straight line, Mather's $\alpha$ function is differentiable at the generator of the line in the first homology.  The next statement gives an algebraic criterion that implies the uniqueness of the asymptotic homologies of Mather measures in the graph. 

\begin{lemma} \label{totally-irrational}
Let $H: T^{*}\mathbb{T}^n \longrightarrow \bR$ a Tonelli Hamiltonian on the torus and let $W \subset T^{*}\mathbb{T}^{n}$ be a $C^{2}$ Lagrangian graph invariant by the Hamiltonian flow. Suppose that all Mather measures have totally irrational rotation number. Then all Mather measures have the same totally irrational rotation number. In particular, Mather's $\alpha$ function is differentiable at the homology class of the measures. 
\end{lemma}

\begin{proof}

We shall prove the statement for $n=3$ since the argument extends straightforwardly to any dimension. Let $v_{1}, v_{2}\mbox{ } \in H_{1}(\mathbb{T}^3 , \mathbb{R})$ two totally irrational vectors, namely, for every base $e_{1}, e_{2}, e_{3} \in H_{1}(\mathbb{T}^3 , \mathbb{Z})$ we have that 
$$ v_{i} = a_{i}e_{1} + b_{i}e_{2} + c_{i}e_{3} $$ 
were the coefficients $a_{1}, b_{i}, c_{i}$ satisfy the following condition: the only linear combination of the type
$$q_{1,i}a_{i} + q_{2,i}b_{i} + q_{3,i} c_{i} = 0$$
for $q_{j, i}$ rational numbers is the trivial one,  $q_{j, i}=0$ for every $j = 1,2$, $i= 1,2,3$. 

Suppose by contradiction, that there exist two Mather measures $\mu_{1}$, $\mu_{2}$ whose homology classes are $v_{1}$, $v_{2}$ respectively, and that $v_{1}$, $v_{2}$ are linearly independent. This means that the matrix 
$$ \begin{pmatrix} 
			a_{1} & a_{2} \\
			b_{1} & b_{2} \\
c_{1} &  c_{2}
			\end{pmatrix}
$$
has rank two. 

Suppose for instance that the matrix $ \begin{pmatrix} 
			a_{1} & a_{2} \\
			b_{1} & b_{2} 
			\end{pmatrix}$ has nonzero determinant. Let $Q_{1},Q_{2}$ be a pair of rational numbers, let us show that there exist $\alpha, \beta$ such that 
$$ \alpha v_{1} + \beta v_{2} = Q_{1}e_{1} + Q_{2}e_{2} + Q_{3}e_{3}$$ 
for some $Q_{3} \in \mathbb{R}$. To see this, we solve the linear system 
$$ \begin{pmatrix} 
			a_{1} & a_{2} \\
			b_{1} & b_{2} \\
c_{1} &  c_{2}
			\end{pmatrix}
\begin{pmatrix} 
			\alpha \\
			\beta 
			\end{pmatrix} = \begin{pmatrix} 
			Q_{1} \\
			Q_{2} \\
Q_{3}
			\end{pmatrix}
$$

Since the matrix $ \begin{pmatrix} 
			a_{1} & a_{2} \\
			b_{1} & b_{2} 
			\end{pmatrix} $ is invertible, the system $$ \begin{pmatrix} 
			a_{1} & a_{2} \\
			b_{1} & b_{2} 
			\end{pmatrix} \begin{pmatrix} 
			\alpha \\
			\beta 
			\end{pmatrix} = \begin{pmatrix} 
			Q_{1} \\
			Q_{2} 
			\end{pmatrix}
$$
has a unique solution. This solution gives the value of $Q_{3}$ in terms of $Q_{1}, Q_{2}, \alpha, \beta$. 

Now, notice that the vector $w= Q_{1}e_{1}+ Q_{2} e_{2} + Q_{3}e_{3}$ is not totally irrational, nor any of its multiples. Therefore, the measure 
$$ \mu = \frac{1}{(\alpha + \beta)}(\alpha \mu_{1} + \beta \mu_{2} )$$ 
is a Mather measure with rotation vector $\frac{w}{(\alpha + \beta)}$, so its rotation vector is not totally irrational, contradicting the assumption of the lemma. Clearly, the same argument extends to any of the $2\times2$ minors of $ \begin{pmatrix} 
			a_{1} & a_{2} \\
			b_{1} & b_{2} \\
c_{1} &  c_{2}
			\end{pmatrix}$ with nonzero determinant. 

The differentiability of Mather's $\alpha$ function in the level set of the Lagrangian submanifold follows from the fundamental properties of the function.

\end{proof}

\section{On the non-existence of Hedlund Lagrangian invariant tori}\label{secHedlund}

Let us recall the notion of {\it Hedlund Lagrangian invariant torus}, introduced in \cite{kn:CarRug}. 

\begin{definition}\label{Hedlundtorus}
A $C^{1}$ Lagrangian torus $S\subset T^*\bT^n$ that is invariant by the Hamiltonian flow of a Tonelli Hamiltonian 
$H : T^*\mathbb{T}^{n} \longrightarrow \mathbb{R}$ is called a {\rm Hedlund Lagrangian torus} if $S$ posseses exactly $n$ closed minimizing orbits whose canonical projections generate the first homology group of $\mathbb{T}^{n}$. 
\end{definition}

This definition is inspired in a famous example due to Hedlund \cite{kn:Hedlund}  of a Riemannian metric on the $n$-torus, with $n \geq 3$, with exactly $n$ closed minimizing geodesics. It was conjectured in \cite{kn:CarRug} that Hedlund Lagrangian invariant tori do not exist in supercritical energy levels; one of the main results in \cite[Proposition 4.1]{kn:CarRug} is the following:

\begin{proposition} \label{Hedlund-tori-nonex1}
Let $W$ be a $C^{1}$ Lagrangian torus invariant by the geodesic flow of a Riemannian metric in $\mathbb{T}^{n}$ containing $n$  orbits whose canonical projections have well defined, linearly independent asymptotic homologies. Then, $W$ is a graph with respect to the canonical projection. 
\end{proposition}

Combining this statement with the generic non-existence of Lagrangian invariant graphs in the $C^{1}$ topology \cite{kn:Ruggiero-no-graphs} we conclude that Hedlund Lagrangian tori do not exist generically in the $C^{1}$ topology for geodesic flows. \\

In the light of Theorem \ref{main2}, we can deduce the non-existence of  Hedlund Lagrangian tori, under the assumptions of Theorem \ref{main2}, without the need of any genericity assumption. Indeed, by Corollary \ref{fibred-torus} and the arguments in Subsection \ref{secconehom} (Case II in the proof of Theorem \ref{cone-fib}), we know that the restriction of the geodesic flow to a $n$-dimensional Lagrangian graph $W$ is conjugate to the suspension of a diffeomorphism $f_{W} : \mathbb{T}^{n-1} \longrightarrow \mathbb{T}^{n-1}$. Let us denote by $\mathcal{W}_{n,I}$ the collection of $C^{1}$ Lagragian graphs invariant by the geodesic flow of a Riemannian metric in $\mathbb{T}^{n}$ such that the diffeomorphism $f_{W}$ is homotopic to the identity for every $W \in \mathcal{W}_{n,I}$. The main result of the subsection is the following:

\begin{proposition} \label{propmain3}
Let $n=3$, then the set of Lagragian graphs in $\mathcal{W}_{n,I}$ such that the collection of Mather measures just contains three elements with linearly independent asymptotic homologies is empty. In particular, there are no Hedlund Lagrangian tori in $\mathcal{W}_{3,I}$. 
\end{proposition}

\begin{proof}

The ideas of the proof are similar to the ideas of the proof of Theorem 5.1. Let $W \in \mathcal{W}_{3,I}$ and suppose that there exist three orbits in $W$, $O_{1}$, $O_{2}$, $O_{3}$ with well defined, independent asymptotic homologies. Let $o_{1}$, $o_{2}$, $o_{3}$ be the corresponding orbits of the diffeomorphism $f_{W}$. 

We claim that the rotation vectors of the orbits $o_{i}$ are in general position in the first homology group $H_{1}(\mathbb{T}^{2}, \mathbb{R})$. Otherwise, the rotation vectors of the 3 orbits would lie in a line in the homology group, and by Lemma \ref{product-hom} the homologies of the orbits $O_{1}$, $O_{2}$, $O_{3}$ would be contained in a plane in the first homology group $H_{1}(\mathbb{T}^{3}, \mathbb{R})$, that is clearly a contradiction since the asymptotic homologies of the orbits  $O_{1}$, $O_{2}$, $O_{3}$ generate $H_{1}(\mathbb{T}^{3}, \mathbb{R})$. 

Therefore, under the assumptions of the proposition, the rotation set of $f_{W}$ contains a convex open set of vectors, and by Franks' work each interior vector in the rotation set gives rise to a periodic orbit of $f_{W}$. Such periodic orbits give rise to closed geodesics in $W$ that are global minimizers since $W$ is a graph. Each closed orbit is then the support of a Mather measure. But by assumption, the collection of Mather measures of $W$ has just 3 elements, yielding a contradiction.  

\end{proof}

\appendix

\section{Mather's theory for Tonelli Lagrangian systems}\label{sec2.4}

In this appendix  we describe Mather's theory  for general Tonelli Lagrangians on compact manifolds. As we have already said before, we refer the reader to \cite{Gonzalonotes, Sorrentinobook} for all the proofs and for a more detailed presentation of this theory.\\

Let $M$ be a compact and connected smooth manifold without boundary.
Denote by $TM$ its tangent bundle and $T^*M$ the cotangent one. A
point of $TM$ will be denoted by $(x,v)$, where $x\in M$ and $v\in
T_xM$, and a point of $T^*M$ by $(x,p)$, where $p\in T_x^*M$ is a
linear form on the vector space $T_xM$. Let us fix a Riemannian
metric $g$ on it and denote by $d$ the induced metric on $M$; let
$\|\cdot\|_x$ be the norm induced by $g$ on $T_xM$; we shall use the same notation for the norm induced on
$T_x^*M$.\\

\begin{definition}\label{defTonelliLag}
A function $L:\,TM\, \longrightarrow \,\bR$ is called a {\em Tonelli
Lagrangian} if:
\begin{itemize}
\item[i)]   $L\in C^2(TM)$;
\item[ii)]  $L$ is strictly convex in the fibres, in the $C^2$ sense, {\it i.e.}, the second partial vertical derivative
${\partial^2 L}/{\partial v^2}(x,v)$ is positive definite, as a quadratic form, for all $(x,v)$;
\item[iii)] $L$ is superlinear in each fibre, {\it i.e.},
            $$\lim_{\|v\|_x\rightarrow +\infty} \frac{L(x,v)}{\|v\|_x} = + \infty.$$
            This condition is equivalent to ask that for each $A\in \bR$ there exists $B(A)\in\bR$ such that
            $$ L(x,v) \geq A\|v\| - B(A) \qquad \forall\,(x,v)\in T M\,.$$ 
\end{itemize}
\end{definition}
Observe that since the manifold is compact, then condition {\it iii)} is independent of the choice of the Riemannian metric $g$.\\

\begin{remark}
One can consider the {\it Hamiltonian} associated to $L$, which is defined  on the cotangent bundle $T^* M$  as its Fenchel transform (or {\it Fenchel-Legendre transform}), see \cite{Rockafellar}: 
\begin{eqnarray} \label{FenchelLeg} H:\; T^*M &\longrightarrow & \bR \nonumber \\
(x,p) &\longmapsto & \sup_{v\in T_xM} \{\langle p,\,v \rangle_x -
L(x,v)\}\, \end{eqnarray} where $\langle \,\cdot,\,\cdot\, \rangle_x$
denotes the canonical pairing between the tangent and cotangent
bundles.
One can easily prove that $H$ is finite everywhere (as a consequence of the superlinearity of $L$), superlinear and strictly convex in each fibre (in the $C^2$ sense). 
Observe that $H$ is also $C^2$. \\
\end{remark}

\medskip
Let $\calM(L)$ be the space of probability measures $\mu$ on 
$TM$  that are invariant  under the Euler-Lagrange flow of $L$ and such that $\int_{TM} \|v\|\,d\m<\infty$.
We will hereafter assume that $\calM(L)$ is endowed with the {\it vague topology}, \ie the weak$^*$--topology induced by the space $C^0_{\ell}$ of continuous functions $f:T M \longrightarrow \R$ having at most linear growth:
$$
\sup_{(x,v)\in T M} \frac{|f(x,v)|}{1+\|v\|} <+\infty\,.
$$

In the case of an autonomous Tonelli Lagrangian, it is easy to see that $\calM(L)$ is non-empty (actually it contains infinitely many measures with distinct supports). In fact,  recall that because of the conservation of the energy
$E(x,v):=H \circ \cL(x,v) = \left\langle \frac{\dpr L }{\dpr v}(x,v),\,v \right\rangle_x - L(x,v)$ along the orbits, each energy level of $E$ is compact (it follows from  the superlinearity condition) and invariant under $\Phi^L_t$. It is a classical result in ergodic theory (sometimes called Kryloff--Bogoliouboff theorem) that a flow on a compact metric space has at least an invariant probability measure, which belongs indeed to $\calM(L)$.

To each  $\m \in \calM(L)$, we may  associate  its {\it average action}:
$$ A_L(\m) = \int_{TM} L\,d\m\,. $$

The action functional 
$A_L : \calM(L) \longrightarrow \R$ is lower semicontinuous with the vague topology on $\calM(L)$.
In particular, this implies that there exists $\m\in\calM(L)$, which minimizes $A_L$ over $\calM(L)$.

\begin{definition}
A measure $\mu\in\calM(L)$, such that $A_L(\mu)=\min_{\calM(L)}A_L$, is called an {\it action-minimizing measure} of $L$.
\end{definition}

Following an idea by John Mather \cite{Ma},  by modifying the Lagrangian (without changing the Euler-Lagrange flow) one can find many other interesting measures besides those found by minimizing $A_L$. 
A similar idea can be implemented for a general Tonelli Lagrangian.
Observe, in fact,  that if $\eta$ is a $1$-form on $M$, we can interpret it as a function 
on the tangent bundle (linear on each fiber)
\begin{eqnarray*} 
\hat{\eta}: TM &\longrightarrow&  \R \\
(x,v) &\longmapsto& \langle \eta(x),\, v\rangle_x
\end{eqnarray*}
and consider a new Tonelli Lagrangian $L_{\eta}:= L - \hat{\eta}$. The 
associated Hamiltonian will be given by
$H_{\eta}(x,p) = H(x,\eta(x) + p)$.

Observe that:
\begin{itemize}
\item[i)] If $\eta$ is closed, then $L$ and $L_{\eta}$ 
have the same Euler-Lagrange flow on $T M$. See \cite{Ma}. 
\item[ii)] If $\m \in \calM(L)$ and $\eta=df$ is an exact $1$-form, then 
$\int{\widehat{df}} d\m =0$. Thus, for a fixed $L$, the minimizing measures will 
depend only on the de Rham cohomology class $c=[\eta] \in H^1(M;\R)$. 
\end{itemize}

Therefore, instead of studying the action minimizing properties of a single Lagrangian, one can consider a family of such ``modified'' Lagrangians, parameterized over $H^1(M;\R)$. 
Hereafter, for any given $c\in H^1(M;\R)$, we will denote by $\eta_c$ a closed $1$-form with that cohomology class.\\

\begin{definition}
Let $\eta_c$ be a closed $1$-form of cohomology class $c$. Then, if $\m \in \calM(L)$ minimizes 
$A_{L_{\eta_c}}$ over $\calM(L)$, we will say that $\m$ is a {\it $c$-action minimizing measure} (or $c$-{\it minimal measure}, or {\it Mather measure} with cohomology $c$).
\end{definition}

One can consider the following function on $H^1(M;\R)$ (the {minus} sign is introduced for a convention that will probably become clearer later on):
\begin{eqnarray}\label{defalpha}
\a: H^1(M;\R) &\longrightarrow& \R \nonumber\\
c &\longmapsto& - \min_{\m\in\calM(L)} A_{L_{\eta_c}}(\m)\,.
\end{eqnarray}
This function $\alpha$ is well-defined (it does not depend on the choice of the representatives of the cohomology classes) and it is easy to see that it is convex. This is generally known as {\it Mather's $\alpha$-function}.
 In particular, it can be proven that  $\a(c)$ is related to the energy level containing such $c$-action minimizing measures  \cite{Carneiro}.\\

We will denote by $\calM_c(L)$ the subset of $c$-action minimizing measures:
$$\calM_c := \calM_c(L)= \{\m \in \calM(L): \;  A_{L_{\eta_c}}(\m)
=-\a (c)\}.$$

We can now define a first important family  of invariant sets: 
the {\it Mather sets}. 

\begin{definition}
For a cohomology class $c \in H^1(M;\R)$, we define 
the {\it Mather set of cohomology class} $c$  as:
\begin{equation} \widetilde{\cM}_c := {\bigcup_{\m \in \calM_c} {\rm supp}\,\m} 
\subset TM\,.\label{2.3}
\end{equation}
The projection on the base manifold $\cM_c = \pi \left(\widetilde{\cM}_c\right)
\subseteq M$ is called {\it projected Mather set} (with cohomology class $c$).
\end{definition}

Properties of this set:
\begin{itemize}
\item[i)] It is non-empty, compact and invariant \cite{Ma}.
\item[ii)] It is contained in the energy level corresponding to $\alpha( c)$ \cite{Carneiro}.
\item [iii)] In \cite{Ma} Mather proved the celebrated {\it graph theorem}: \\
{\it Let $\pi: T M \longrightarrow M$ denote the canonical projection. Then,
 $\pi|{\widetilde{\cM}_c}$ is an injective 
mapping of $\widetilde{\cM}_c$ into $M$, and its inverse $\pi^{-1}: \cM_c 
\longrightarrow \widetilde{\cM}_c$ is 
Lipschitz.}\\
 \end{itemize}

\vspace{10 pt}
Now, we would like to shift our attention to a related problem. Instead of considering different minimizing problems over $\calM(L)$, obtained by modifying the Lagrangian $L$, one can alternatively try to minimize the Lagrangian $L$ by putting some {constraint}, such as, for instance, fixing the {\it rotation vector} of the measures. In order to generalize this to Tonelli Lagrangians on compact manifolds, we first need to define what we mean by rotation vector of an invariant measure.

Let $\mu\in \calM(L)$. Thanks to the superlinearity of $L$, 
the integral  $ \int_{TM} \hat{\eta} d\m$ 
is well defined and finite for any 
closed 1-form $\eta$ on $M$.
Moreover,  if $\eta$ is exact, then this integral is zero, \ie 
$ \int_{TM} \hat{\eta} d\m=0$.
Therefore, one can define a linear functional: 
\begin{eqnarray*}
H^1(M;\R) &\longrightarrow& \R \\
c &\longmapsto& \int_{TM} \hat{\eta} d\m\,,
\end{eqnarray*}
where $\eta$ is any closed $1$-form on $M$ with cohomology class $c$. By 
duality, there 
exists $\rho (\m) \in H_1(M;\R)$ such that
$$
\int_{TM} \hat{\eta} \,d\m = \langle c,\rho(\m) \rangle
\qquad \forall\,c\in H^1(M;\R)$$ 
(the bracket on the right--hand side denotes the canonical pairing between 
cohomology and 
homology). We call $\rho(\m)$ the {\it rotation vector} of $\m$. This rotation vector is 
the same as the Schwartzman's asymptotic cycle of $\mu$ (see \cite{kn:Schwartz} and \cite{Sorrentinobook} for more details).

Using that the action functional $A_L: \calM(L) 
\longrightarrow \R$ is lower semicontinuous, one can prove that the map $\rho: \calM(L) \longrightarrow H_1(M;\R)$ is continuous and surjective, 
{\it i.e.}, for every  $h\in H_1(M;\R)$  there exists $\m\in \calM(L)$ with  $A_L(\m) < \infty$ and $\rho(\m)=h$ (see \cite{Ma}).\\

Following Mather \cite{Ma}, let us consider the minimal value of the average action $A_L$ over the 
probability measures with rotation vector $h$. Observe that this minimum is actually achieved because of the lower semicontinuity of $A_L$ and the compactness of $\rho^{-1}(h)$ ($\rho$ is continuous and $L$ superlinear). Let us define
\begin{eqnarray} \label{defbeta}
\b: H_1(M;\R) &\longrightarrow& \R  \nonumber\\
h &\longmapsto& \min_{\m\in\calM(L):\,\rho(\m)=h} A_L(\m)\,.\label{2.2}
\end{eqnarray}
This function $\beta$ is what is generally known as {\it Mather's 
$\beta$-function} and it is immediate to check that it is convex; this function is sometime called {\it effective Lagrangian}.\\

We can now define what we mean by action minimizing measure with a given rotation vector.

\begin{definition}
A measure $\m \in \calM(L)$ realizing the minimum in (\ref{2.2}), \ie such that $A_L(\m)~=~\b(\rho(\m))$, is called an {\it action minimizing} (or {\it minimal}, or {\it Mather}) {\it measure} with rotation vector $\rho(\m)$.
\end{definition}

 We will denote by $ \calM^h(L)$ the subset of action minimizing measures with rotation vector $h$:
{\small $$ \calM^h := \calM^h(L) = \{\m \in \calM(L):  \;\rho(\m)= h \; 
{\rm and}\;  A_L(\m)=\beta(h)\}.$$}

 This allows us to define another important familty of invariant sets.

 \begin{definition}
 For a homology class (or rotation vector) $h\in H_1(M;\R)$, we define the
{\it Mather set corresponding to a rotation vector} $h$ as
\begin{equation} \widetilde{\cM}^h := {\bigcup_{\m \in \calM^h} {\rm supp}\,\m} 
\subset TM\,,\label{2.5}\end{equation}
and the projected one as $\cM^h = \pi \left(\widetilde{\cM}^h\right) 
\subseteq M$. 
 \end{definition}

 Similarly to what we have already seen above, this set satisfies the following properties:
\begin{itemize}
\item[i)] It is non-empty, compact and invariant.
\item[ii)] It is contained in a given energy level.
\item [iii)] It also satisfies the {\it graph theorem}:\\ {\it let $\pi: T M \longrightarrow M$ denote the canonical projection. Then,
 $\pi|{\widetilde{\cM}^h}$ is an injective 
mapping of $\widetilde{\cM}^h$ into $M$, and its inverse $\pi^{-1}: \cM^h 
\longrightarrow \widetilde{\cM}^h$ is 
Lipschitz.}\\
 \end{itemize}
 
 \vspace{10 pt}

\begin{remark}\label{rem2.17}
({\it i}) In the above discussion we have only discussed properties of invariant probability measures associated to the system. Actually, one could study directly orbits of the systems and look for orbits that {globally minimize} the action of a modified Lagrangian (in the same spirit as before). This would lead to the definition of two other families of invariant compact sets, the {\it Aubry sets} $\widetilde{\cA}_c$  and the {\it Ma\~n\'e sets} $\widetilde{\cN}_c$, which are also parameterized by $H^1(M;\R)$ (the parameter which describes the modification of the Lagrangian, exactly in the same way  as before). For a given $c\in H^1(M;\R)$, these sets contain  the Mather set $\widetilde{\cM}_c$, and this inclusion may be strict. In fact, while the motion on the Mather sets is {\it recurrent} (it is the union of the supports of invariant probability measures), the Aubry and the Ma\~n\'e sets may contain non-recurrent orbits as well.\\
({\it ii }) Differently from what happens with invariant probability measures, it will not be always possible to find {\it action-minimizing orbits} for any given rotation vector (not even possible to define a rotation vector for every action minimizing orbit). For instance, an example due to Hedlund \cite{kn:Hedlund} provides the existence of a Riemannian metric on a three-dimensional torus, for which minimal geodesics exist only in three directions. The same construction can be extended to any dimension larger than three.\\
\end{remark}

The above discussion led to two equivalent formulations of the minimality 
of an invariant probability measure $\m$:
\begin{itemize}
\item there exists a homology class $h \in {\rm H}_1(M;\R)$, namely its 
rotation vector $\rho(\m)$, such that $\m$ minimizes $A_L$ amongst all 
measures in $\calM(L)$ with rotation vector $h$, \ie $A_L(\m)=\b (h)$.
\item There exists a cohomology class $c \in {\rm H}^1(M;\R)$, such that $\m$ minimizes $A_{L_{\eta_c}}$ 
amongst all probability measures in $\calM(L)$, \ie $A_{L_{\eta_c}}(\m)=-\a (c)$.\\
\end{itemize}

\noindent What is the relation between two these different approaches? Are they equivalent, \ie 
$\bigcup_{h \in {\rm H}_1(M;\R)} \calM^h = \bigcup_{c \in {\rm H}^1(M;\R)} 
\calM_c\,$ ?\\

In order to comprehend the relation between these two families of action-minimizing measures, we need to understand better the properties of the these two functions that we have introduced above:
$$\a: H^1(M;\R) \longrightarrow \R \quad {\rm and}\quad  \b:H_1(M;\R) \longrightarrow \R.$$

\begin{remark} Before stating  the main relation between these two functions, let us recall some definitions and results from classical convex analysis (see \cite{Rockafellar}).\\
{\bf (i)}  Given a convex function $\varphi: V \longrightarrow \R\cup \{+\infty\}$ on a finite dimensional vector space $V$, one can consider a {\it dual} (or {\it conjugate}) function  defined on the dual space $V^*$, via the so-called {\it Fenchel transform}: $\varphi^*(p):= \sup_{v\in V} \big(p\cdot v - \varphi(v)\big)$.  In our case, the following holds.\\
{\bf (ii)} Like any convex 
function on a
finite-dimensional space, $\b$
admits a subderivative at each point $h\in H_1(M;\R)$, \ie we can find $c\in 
H^1(M;\R)$ such that
\begin{eqnarray*}
\forall h'\in H_1(M;\R), \quad \b(h')-\b(h)\geq \langle c,h'-h\rangle.\end{eqnarray*}
We will denote by $\partial \b(h)$ the set of $c\in 
H^1(M;\R)$ that 
are subderivatives of $\b$ at $h$, \ie the set of $c$'s which satisfy the above  
inequality.  Similarly, we will denote by $\dpr \a(c)$ the set of subderivatives of $\a$ at $c$.
Actually, Fenchel's duality implies an easier characterization of subdifferentials:
{ \it $c\in \partial \b(h)$ if and only if $\langle c,h\rangle=\a(c )+\b(h)$} (similarly for $h\in \partial \a(c)$).\\
\end{remark}

\begin{proposition}\label{alphabetasuperlinear}
$\a$ and $\b$ are convex conjugate, \ie
$\a^* = \b$ and $\b^* = \a$. In particular, it follows that $\a$ and $\b$ have superlinear growth.\\
Moreover,  if  $\m \in \calM(L)$ is an invariant probability measure, then:
\begin{itemize}
\item[{\rm (i)}] $A_L(\m)=\b(\rho(\m))$ {if and only if}  there exists $c\in {\rm H}^1(M;\R)$ such that $\m$ minimizes $A_{L_{\eta_c}}$  
{\rm(}\ie $A_{L_{\eta_c}}(\mu)=-\a(c)${\rm)}.
\item[{\rm (ii)}] If $\m$ satisfies $A_L(\m)=\b(\rho(\m))$ and $c \in 
H^1(M;\R)$, then $\m$ minimizes $A_{L_{\eta_c}}$ if and only if $c\in \partial 
\b(\rho(\mu))$ {\rm(}or equivalently $\langle c,h\rangle=\a(c)+\b(\rho(\m)\rm{)}$.\\
\end{itemize}
\end{proposition}

\begin{remark}\label{remarkinclusionsmathersets}
{\bf (i)} It follows from the above proposition that both minimizing procedures lead to the same sets of invariant probability measures:
$$
\bigcup_{h \in {\rm H}_1(M;\R)} \calM^h = \bigcup_{c \in {\rm H}^1(M;\R)} 
\calM_c\,.
$$
In other words,  minimizing over the set of invariant measures with a fixed rotation vector or globally minimizing the modified Lagrangian (corresponding to a certain cohomology class) are dual problems,  as the ones that often appears in linear programming and optimization. In some sense, modifying the Lagrangian by a closed $1$-form 
is analog to the method of Lagrange multipliers for searching constrained critical points of a function.\\
{\bf (ii)} In particular we have the following inclusions between Mather sets:
$$
 c\in \dpr \beta(h) \quad \Longleftrightarrow \quad h \in \partial \a(c) \quad \Longleftrightarrow \quad \widetilde{\cM}^h \subseteq \widetilde{\cM}_c\,.
 $$
Moreover, for any $c\in H^1(M;\R)$:
$$
\widetilde{\cM}_c = \bigcup_{h\in \dpr \a(c)} \widetilde{\cM}^h\,.
$$
Observe that the non-differentiability of $\alpha$ at some $c$ produces the presence in $ \widetilde{\cM}_c$ of (ergodic) invariant probability measures with different rotation vectors. On the other hand, the non-differentiability of $\beta$ at some $h$ implies that there exist $c\neq c'$ such that $ \widetilde{\cM}_c \cap  \widetilde{\cM}_{c'}\neq \emptyset$).\\
{\bf (iii)} The minimum of the $\a$-function is sometime called {\it Ma\~n\'e's strict critical value}. Observe that if $\a(c_0) = \min \a(c)$, then $0\in \dpr \a(c_0)$ and $\beta(0)=-\a(c_0)$. Therefore, the measures with zero homology are contained in the least possible energy level containing Mather sets and $\widetilde{\cM}^0 \subseteq \widetilde{\cM}_{c_0}$. This inclusion might be strict, unless $\a$ is differentiable at $c_0$; in fact, there may be other action minimizing measures with non-zero rotation vectors corresponding to the other subderivatives of $\a$ at $c_0$.\\

\end{remark}

\vspace{10 pt}

\section{Corrigendum to   \cite[Lemma 3.3]{kn:CarRug}} \label{Appcorrigendum}

In this appendix, we provide a correction  to the proof of  \cite[ Lemma 3.3]{kn:CarRug}, whose statement should read as follows:
\bigskip

\begin{lemma} \label{full-pi2}
Let $M$ be a compact smooth manifold, and let $H : T^{*}M \longrightarrow \mathbb{R}$ be a Tonelli Hamiltonian. Let $\pi : T^{*}M \longrightarrow M$ be the canonical projection, let $W \subset T^{*}M$ be a $C^{2}$ Lagrangian submanifold such that the natural homomorphism 
	$$\pi_{*} : H_{1}(W, \mathbb{Z}) \longrightarrow H_{1}(M, \mathbb{Z})$$ 
induced by $\pi$ is an epimorphism. Then the natural homomorphism induced by $\pi$ between the first cohomology groups 
$$\pi^{*} : H^{1}(M, \mathbb{Z}) \longrightarrow H^{1}(W, \mathbb{Z})$$ 
given by $\pi^{*}([\beta]) = [\pi^{*}(\beta)]$, is a monomorphism. 

In particular, if $M = \mathbb{T}^{n}$ and $W$ is a $C^{2}$ $n$-dimensional torus, then the natural homomorphism $\pi^{*}$ is an isomorphism.
	\end{lemma}

\begin{proof}

Let us start as in \cite{kn:CarRug}, by recalling the duality between the first homology group $H_{1}(N,\mathbb{R})$ and the first cohomology group $H^{1}(N,\mathbb{R})$ given by the Hurewicz homomorphism in a smooth manifold $N$: we can represent cohomology classes as homomorphisms $\rho : \pi_{1}(N) \longrightarrow \mathbb{R}$. Namely, given a closed 1-form $\omega$ defined in $N$, the following map: 
$$ \rho_{\omega} (\alpha) := \int_{\alpha} \omega $$ 
defined in the set of closed curves, induces a homomorphism $\rho_{\omega} : \pi_{1}(N) \longrightarrow \mathbb{R}$ 
in the first homotopy group. By Stokes Theorem, this map is well defined in the first homology group, namely, if $\sigma$ is a closed curve homologous to $\alpha$, then $ \rho_{\omega} (\alpha) = \rho_{\omega} (\sigma) $, so $\rho_{\omega}$ induces as well a homomorphism $\rho_{\omega} : H_{1}(N,\mathbb{R}) \longrightarrow \mathbb{R}$. Therefore, $\rho_{\omega}$ represents a cohomology class of $N$. Let us denote by $\Hom (\pi_{1}(N), \mathbb{R})$ the group of homomorphisms from $\pi_{1}(N)$ to $\mathbb{R}$. 

Now, the canonical projection defines a natural map 
$$ \hat{\pi}_{*} : \Hom(\pi_{1}(M), \mathbb{R}) \longrightarrow \Hom(\pi_{1}(W), \mathbb{R})$$
given by $\hat{\pi}_{*}(\rho_{\bar{\omega}}) = \rho_{\pi^{*}(\bar{\omega})}$, where $\bar{\omega}$ is a closed 1-form 
in $M$. The map $\hat{\pi}_{*}(\rho_{\bar{\omega}})$ evaluated on a closed curve $\alpha  \subset W$ is 
$$ \hat{\pi}_{*}(\rho_{\bar{\omega}})(\alpha) = \rho_{\pi^{*}(\omega)}(\pi (\alpha)) .$$
This map is well defined in the homology class of $\alpha$, so as a map between first homology groups it is just 
$$ \hat{\pi}_{*}(\rho_{\bar{\omega}})([\alpha]) = \rho_{\pi^{*}(\bar{\omega})}(\pi_{*}[ (\alpha)]),$$
where $[\alpha] \in H_{1}(W,\mathbb{Z})$ denotes the homology class of $\alpha$. 

By the assumptions of the Lemma, the map 
$$\pi_{*} : H_{1}(W, \mathbb{Z}) \longrightarrow H_{1}(M, \mathbb{Z})$$ 
is an epimorphism, so the values
$$\rho_{\pi^{*}(\bar{\omega})}(\pi_{*}[ (\alpha_{i})]) = \int \bar{\omega} (\pi_{*}[\alpha_{i}])$$  in a basis 
$[ (\alpha_{i})]$ of $H_{1}(W,\mathbb{Z})$ determine uniquely the form $\bar{\omega}$ by the De Rham cohomology theory.  

Hence, the correspondence $\bar{\omega} \longmapsto \pi^{*}(\bar{\omega})$ is well defined in $H^{1}(M,\mathbb{Z})$, and applying once more De Rham cohomology theory we conclude that the map $\pi^{*}$ is injective as claimed. 

When $M$ and $W$ are diffeomorphic to the $n$-torus, the map $\pi^{*}$ is an isomorphism since the first cohomology groups of $H^{1}(M,\mathbb{Z})$ and $H^{1}(W, \mathbb{Z})$ are isomorphic. 

\end{proof}

\end{document}